\documentclass[11pt]{amsart}
\usepackage{amssymb}
\usepackage{amsfonts}
\usepackage{amsmath}
\usepackage{amscd}
\usepackage[matrix,arrow,graph]{xy}
\input{epsf}

\begin{document}

\title{Perverse Sheaves on Grassmannians}
\author{Tom Braden}
\address{Department of Mathematics and Statistics\\ University of Massachusetts\\ Amherst, MA 01003}
\email{braden@math.umass.edu}
\thanks{This research supported in part by NSF grant DMS 9304580}
\begin{abstract}
We compute the category of perverse sheaves on Hermitian symmetric spaces in types A and D,
constructible with respect to the Schubert stratification.  The calculation
is microlocal, and uses the action of the Borel group to study the
geometry of the conormal variety $\Lambda$.
\end{abstract} 
\maketitle
\tableofcontents

\theoremstyle{plain}
\newtheorem{thm}{Theorem}[subsection]
\newtheorem{lemma}[thm]{Lemma}
\newtheorem{prop}[thm]{Proposition}
\newtheorem{cor}[thm]{Corollary}
\theoremstyle{definition}
\newtheorem*{conj}{Conjecture}
\newtheorem*{defn}{Definition}
\newtheorem*{nota}{Notation}
\theoremstyle{remark}
\newtheorem*{rem}{Remark}

\hyphenation{Mac-Pher-son}
\hyphenation{canon-ical}

\newcommand{\ra}{\rightarrow}
\newcommand{\ral}{\mathop{\rightarrow}\limits}
\newcommand{\la}{\leftarrow}
\newcommand{\lra}{\leftrightarrow}
\newcommand{\wt}{\widetilde}

\def\rank{\mathop{\rm rank}\nolimits}
\def\ker{\mathop{\rm Ker}\nolimits}
\def\coker{\mathop{\rm Coker}\nolimits}
\def\im{\mathop{\rm Im}\nolimits}
\def\Id{{\rm Id}}
\def\id{{\rm id}}
\def\ob{\mathop{\rm ob}}
\def\lcm{\mathop{\rm lcm}\limits}
\def\Hom{\mathop{\rm Hom}}
\def\Ext{\mathop{\rm Ext}\nolimits}
\def\Supp{\mathop{\rm supp}}
\def\Span{\mathop{\rm span}}
\newcommand{\half}{\frac{1}{2}}
\newcommand{\be}{{\bf e}}
\newcommand{\mtop}{{\rm top}}

\newcommand{\C}{{\mathbb C}}
\newcommand{\Z}{{\mathbb Z}}
\newcommand{\R}{{\mathbb R}}
\newcommand{\Hi}{{\mathbb H}}
\newcommand{\Rb}{{R_\lambda}}
\newcommand{\Wb}{{\wt{Z}}}
\newcommand{\lambdap}{{\lambda_0}}
\newcommand{\SO}{{S_O}}
\newcommand{\Rge}{\R_{\ge 0}}
\newcommand{\Zge}{\Z_{\ge 0}}
\newcommand{\PP}{{\mathbb P}}
\newcommand{\Q}{{\mathbb Q}}
\newcommand{\F}{{\mathbb F}}
\newcommand{\HH}{{\mathbb H}}
\newcommand{\N}{{\mathbb N}}
\newcommand{\cal}{\mathcal}
\newcommand{\cS}{{\cal S}}
\newcommand{\cO}{{\cal O}}
\newcommand{\cU}{{\cal U}}
\newcommand{\cC}{{\cal C}}
\newcommand{\cF}{{\cal F}}
\newcommand{\cE}{{\cal E}}
\newcommand{\cL}{{\cal L}}
\newcommand{\cP}{{\cal P}}
\newcommand{\cQ}{{\cal Q}}
\newcommand{\cR}{{\cal R}}
\newcommand{\cD}{{\cal D}}
\newcommand{\cA}{{\cal A}}
\newcommand{\cM}{{\cal M}}
\newcommand{\IC}{{\mathbf {IC^{\textstyle \cdot}}}}
\newcommand{\mb}{\mathbf}

\hyphenation{para-metrizing}

In \cite{BB} Beilinson and Bernstein gave an equivalence of
categories between certain categories of 
perverse sheaves on a flag variety $G/B$ and
certain representations of the Lie algebra $\mathfrak g$.  This correspondence
has most often been used to understand irreducible representations
by studying the corresponding perverse sheaves. 
The complete structure of the category of perverse sheaves has been
computed in a few simple cases (in particular, for projective spaces
with the Schubert stratification) but existing techniques have not allowed
computation of more complicated examples.
In this paper we give a quiver description of the category $\cP_\Lambda(X)$
when $X = G/P$ is a Hermitian symmetric space in type $A$ or $D$
(i.e.\ a Grassmannian in type A, or an isotropic Grassmannian in type D), 
stratified
by the Schubert stratification.  Here $\Lambda \subset T^*X$ is
the conormal variety to the stratification; the category $\cP_\Lambda$
of perverse sheaves with characteristic variety contained in $\Lambda$
is the same as the category of Schubert-constructible perverse sheaves.  

Our strategy is to study perverse sheaves microlocally, that is,
as objects supported on $\Lambda$.  Such a description
exists, at least theoretically, via the theory of regular singularities
$\cE$-modules.  A topological description of this category was
given in \cite{GMV}.  In practice, computing such a category proceeds
from smooth points of $\Lambda$ inwards to deeper singularities.  
A conjecture of Kashiwara says that 
only codimension $0$, $1$, and $2$ pieces of $\Lambda$ should be
necessary in the computation; our methods show as a corollary
that this conjecture
holds for our spaces, although the proof is quite special to our
particular geometry.

There are several pleasant features of our varieties which make
the computation of $\cP_\Lambda(X)$ reasonable, all of which
 fail for full flag varieties $G/B$.
First, the action of Borel group $B$ on the conormal variety $\Lambda$ 
has finitely 
many orbits.  This provides a natural stratification of $\Lambda$ and
allows a simple description of the geometry of the strata
and how they intersect.  The author is unaware of even an
algorithm to decide whether two components of the conormal variety to 
the Schubert stratification of a full flag variety
$G/B$ meet in codimension one.

Second, the action of the Borel group $B$ on $\Lambda$ 
has connected stabilizers, so the fundamental groups 
of the orbits are free abelian groups, generated by $\pi_1(B) = \pi_1(T)$.  
Already 
for the full flag variety for $SL_4$ there are smooth components of $\Lambda$
with nonabelian $\pi_1$.  
If $X$ is the the type $B$ Hermitian symmetric space (the Lagrangian 
Grassmannian), the stabilizers are not connected, 
so the fundamental groups of the orbits cannot
be completely described by the action of $\pi_1(T)$.  This is essentially the
reason we do not consider this case in this paper.

Finally, all the singularities of the Schubert stratification of $G/P$ are
conical.  This simplifies things considerably; it means that all the
codimension one intersections of components of $\Lambda$ look like
the conormal variety to a line bundle $L$ stratified as 
$Z \cup (L \setminus Z)$
where $Z$ is the zero section.  It also allows the use of the Fourier transform
to identify microlocal perverse sheaves on different spaces, rather than
 contact transformations, which are harder to compute with.  

A general description of microlocal perverse sheaves 
on the union of codimension zero and one 
strata of $\Lambda$ for conical stratifications was given in
\cite{BG}; in that paper it was applied to stratifications where there
are no codimension two strata.  There are codimension two orbits
for the spaces we consider, but the geometry of $\Lambda$ near
these orbits is as simple as possible.  It is just the conormal
variety to a direct sum of two line bundles $L_1 \oplus L_2$,
with the ``normal crossings" stratification.

The resulting presentation by generators and relations, while directly 
arising from the conormal geometry, is not algebraically the most pleasing.
Mikhail Khovanov has described a quiver algebra \cite{Kh} 
arising from an algebra of cobordisms, which he uses to 
``categorify'' invariants of links and tangles.  His 
algebra is isomorphic to a subquotient of our algebra for the type A
Grassmannian; a proof will appear in \cite{BK}.  Khovanov's algebra is 
naturally graded, 
with quadratic relations, and thus our algebra is also, at least
in the type A case.  Such a grading is a crucial ingredient in 
Koszul duality theory \cite{BGS}.

The paper is organized as follows.  \S1 presents the quiver
category that describes the category ${\cal P}_\Lambda(X)$, and 
describes the simple objects in this category.  After some combinatorial
preliminaries the quiver categories are described in sections 1.4 (type A)
and 1.6 (type D).  \S2 describes
 the geometry and combinatorics of the $B$-orbits of
$\Lambda$.  \S3 introduces microlocal perverse sheaves, describes
the building blocks (monodromic and normal crossings perverse sheaves) 
which are ``glued
together" to give the final answer, and finally \S4 gives the necessary
identifications to carry out the gluing.

\section{Preliminaries and statement of results}

\subsection{The classical (type A) Grassmannian}
Consider the complex Grassmannian $X = X_{k, l}$ parametrizing
$k$-dimensional sub-vector spaces of $\C^{n}$, where $n = k+l$.  Letting 
$G = SL(n, \C)$, then $G$ acts on $X$ transitively, and 
$X = G/P$ for $P = P_{k,l}$ a maximal parabolic subgroup of $G$.
The Borel subgroup $B$ of upper triangular matrices in $G$ acts on
$X$ with finitely many orbits, called Schubert cells.  These cells
give a stratification of $X$ which is described as follows.

Let $\Omega$ 
denote the collection of all partitions, i.e.\ nonincreasing sequences
$\lambda_1 \ge \lambda_2 \ge \dots \ge \lambda_i\ge \dots$ 
of nonnegative integers which are eventually zero.
  Let $\Omega_{k,l} \subset \Omega$
be the set of $\lambda$ for which $\lambda_{k+1} = 0$ and $\lambda_i\le l$ 
for all $i$.  In other words, $\lambda \in \Omega_{k,l}$ if and only
if the Young diagram $\Delta(\lambda)$ of $\lambda$ 
fits in a rectangle with $k$ rows and $l$ columns.  
Here we put 
\[\Delta(\lambda) = \{(i, j) \in \N\times \N \mid j \le \lambda_i\}\]
(we are using matrix coordinates for 
$\N\times \N$: $(i,j)$ is the point in the $i$th row and $j$th column
from the upper left).

Denote the standard flag fixed by $B$ by 
$\C^1 \subset \C^2 \subset \dots \subset \C^n$, and let $\be_1, \dots, \be_n$
be the standard basis of $\C^n$.
For $\lambda \in \Omega_{k,l}$
define the Schubert cell $X_\lambda$ to be
\[X_\lambda = \{\, V \in X_{k,l} \mid \dim(V \cap \C^{d_j}) = j,
 \; 1\le j \le k\},\]
where we put $d_j$ = $\lambda_{k-j+1} + j$.
It is a smooth affine variety isomorphic to $\C^{|\lambda|}$, where
$|\lambda| = \sum \lambda_i$.  Put a partial order on $\Omega$ by
defining $\lambda \le \lambda'$ if and only if $\lambda_i \le \lambda'_i$
for all $i$ --
i.e. by inclusion of Young diagrams.  Then $X_\lambda \subset
\overline{X_{\lambda'}}$ if and only if $\lambda \le \lambda'$.

There is another way of describing a partition $\lambda$ which
will be useful.  Let $\Hi = \Z + \half$ be the set of ``half integers".  
Given $\lambda$, define a function $\phi = \phi_\lambda
\colon \Hi \to \{+1, -1\}$ by letting $\phi(\alpha) = -1$ for
$\alpha \in \{\, \lambda_i - i + \half \mid i>0\,\}$, and 
$\phi(\alpha)=+1$ otherwise.
Geometrically, the function $\phi$ is obtained by moving along the outer
edge of the Young diagram $\Delta(\lambda)$ and giving a $-1$ for every 
step up and a $+1$ for every step to the right.  More precisely, 
$\phi(\alpha)$ is determined by the orientation of the 
boundary segment that intersects the line $y+x = \alpha$.

\begin{figure}[tb]
\begin{center}
\leavevmode
\hbox{
\epsfxsize=3in
\epsffile{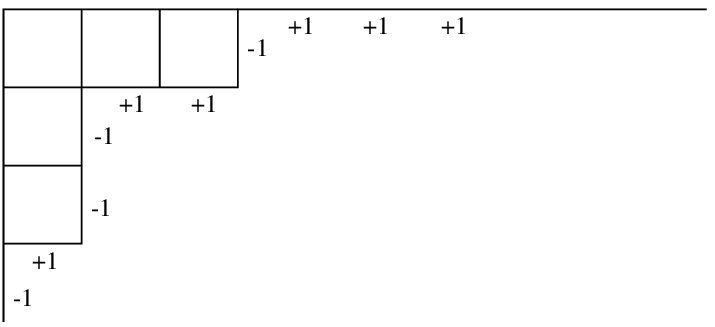}}
\end{center}
\caption{ Computing $\phi_{\lambda}$ by following the boundary of 
$\Delta(\lambda)$}
\label{lbdaphi}
\end{figure}

This produces a 
bijection between $\Omega$ and the set of all functions 
$\phi\colon \Hi \to \{+1, -1\}$ which are $-1$ for all sufficiently negative
integers and $+1$ for all sufficiently positive ones, and for which 
the sum $\sum^\beta_{\alpha=-\beta} \phi(\alpha)$ is zero for all
$\beta$ large enough (See figure \ref{lbdaphi}).

\subsection{$\lambda$-pairs: type A} \label{lprsec}
Let $\wt{\Pi}(\lambda) = 
\phi^{-1}(-1) \times \phi^{-1}(+1) \subset \Hi \times \Hi$.

\begin{defn}
Call a pair $(\alpha, \beta) \in \wt{\Pi}(\lambda)$
a ``$\lambda$-pair'' if $\beta$ is the smallest number 
for which $\beta>\alpha$ 
and $\sum_{\alpha \le \gamma \le \beta} \phi(\gamma) = 0.$
Let $\Pi(\lambda) \subset \wt{\Pi}(\lambda)$ be the set
of $\lambda$-pairs.
\end{defn}

\begin{figure}[tb]
\begin{center}

\setlength{\unitlength}{3947sp}%
\newcommand{\sss}{\scriptscriptstyle}
\begingroup\makeatletter\ifx\SetFigFont\undefined%
\gdef\SetFigFont#1#2#3#4#5{%
  \reset@font\fontsize{#1}{#2pt}%
  \fontfamily{#3}\fontseries{#4}\fontshape{#5}%
  \selectfont}%
\fi\endgroup%
\begin{picture}(3537,1677)(1051,-1348)
\put(1726,-136){$\sss\frac{1}{2}$}
\put(2101,-136){$\sss\frac{3}{2}$}
\put(2220,-50){$\sss\frac{5}{2}$} 
\put(2476,239){$\sss\frac{7}{2}$} 
\put(1241,-886){$\sss -\frac{5}{2}$}
\put(951,-1201){$\sss -\frac{7}{2}$}
\put(1379,-766){$\sss -\frac{3}{2}$}
\put(1379,-421){$\sss -\frac{1}{2}$} 
\thinlines
\put(1201,-1336){\line( 0, 1){1500}}
\put(1201,164){\line( 1, 0){1125}}
\put(1201,-586){\line( 1, 0){375}}
\put(1201,-211){\line( 1, 0){375}}
\put(1576,-211){\line( 0, 1){375}}
\put(1951,-211){\line( 0, 1){375}}
\put(1201,-961){\line( 1, 0){375}}
\put(1576,-961){\line( 0, 1){750}}
\put(1576,-211){\line( 1, 0){750}}
\put(2326,-211){\line( 0, 1){375}}
\put(2326,164){\line( 1, 0){2250}}
\put(1611,-421){\line( 1, 1){170}}
\put(1621,-781){\line( 1, 1){525}}
\put(1246,-1171){\line( 1, 1){170}}
\put(2351,-51){\line( 1, 1){170}}
\end{picture}


\end{center}
\caption{Finding $\Pi([3,1,1])$} 
\label{lampr}
\end{figure}
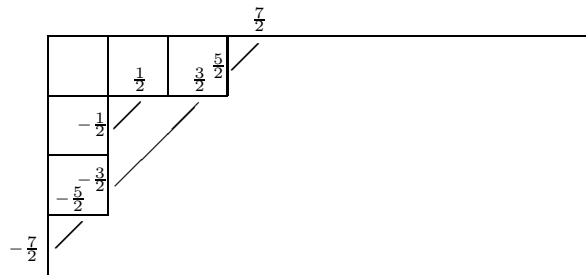

If $(\alpha, \beta)\in \Pi(\lambda)$, and we start from the center of the 
boundary segment of $\Delta(\lambda)$ corresponding to $\alpha$, then
$\beta$ is the first boundary segment encountered by a ray extended 
up and to the right with slope $+1$.  Figure \ref{lampr}
shows that for $\lambda = (3, 1, 1)$, $(-\half, \half), 
(-\frac{3}{2}, \frac{3}{2}),
(-\frac{7}{2}, -\frac{5}{2}),$ and 
$(\frac{5}{2}, \frac{7}{2})$ are all $\lambda$-pairs.

\begin{lemma}
\label{lpnest}  Every $\alpha \in \Hi$ appears in exactly one $\lambda$-pair.
The $\lambda$-pairs define a set of nested intervals; if $(\alpha, \beta)$
and $(\alpha', \beta')$ are $\lambda$-pairs, and $\alpha' \in [\alpha, \beta]$,
then $\beta' \in [\alpha, \beta]$ also. 
\end{lemma}
Put a partial order on $\wt{\Pi}(\lambda)$ by letting $(\alpha, \beta) \le
(\alpha', \beta')$ whenever $\alpha' \le \alpha$ and $\beta \le \beta'$, and
let $\Pi(\lambda)$ have the induced partial order.
 Lemma $\ref{lpnest}$ implies that
 $\Pi(\lambda)$ is a tree; each 
element has a unique smallest element dominating it (its ``parent").

Given a $\lambda$-pair $\pi= (\alpha, \beta)$, define a new partition 
$\lambda'$ by letting 
$\phi_{\lambda'}(\alpha) = +1$, $\phi_{\lambda'}(\beta) = -1$, and
$\phi_{\lambda'}(\gamma) = \phi_\lambda(\gamma)$ otherwise (see Figure 
\ref{lplp}).
We will denote this relation by $\lambda \ral^\pi \lambda'$;
if it holds for some $\lambda$-pair, we write $\lambda \ra \lambda'$, 
and if $\lambda \ra \lambda'$ or $\lambda' \ra \lambda$, we write
$\lambda \lra \lambda'$.

\begin{figure}[tbh] 
\setlength{\unitlength}{3947sp}%
\begingroup\makeatletter\ifx\SetFigFont\undefined%
\gdef\SetFigFont#1#2#3#4#5{%
  \reset@font\fontsize{#1}{#2pt}%
  \fontfamily{#3}\fontseries{#4}\fontshape{#5}%
  \selectfont}%
\fi\endgroup%
\begin{picture}(3612,1524)(976,-1348)
\thinlines
\put(1201,-1336){\line( 0, 1){1500}}
\put(1201,164){\line( 1, 0){1125}}
\put(1201,-586){\line( 1, 0){375}}
\put(1201,-211){\line( 1, 0){375}}
\put(1576,-211){\line( 0, 1){375}}
\put(1951,-211){\line( 0, 1){375}}
\put(1201,-961){\line( 1, 0){375}}
\put(1576,-961){\line( 0, 1){750}}
\put(1576,-211){\line( 1, 0){750}}
\put(2326,-211){\line( 0, 1){375}}
\put(2326,164){\line( 1, 0){2250}}
\multiput(1576,-961)(107.14286,0.00000){4}{\line( 1, 0){ 53.571}}
\multiput(1951,-961)(0.00000,107.14286){4}{\line( 0, 1){ 53.571}}
\multiput(1951,-586)(107.14286,0.00000){4}{\line( 1, 0){ 53.571}}
\multiput(2326,-586)(0.00000,107.14286){4}{\line( 0, 1){ 53.571}}
\multiput(1576,-586)(107.14286,0.00000){4}{\line( 1, 0){ 53.571}}
\multiput(1951,-586)(0.00000,107.14286){4}{\line( 0, 1){ 53.571}}
\put(2161,-911){\Large$\lambda'$}
\put(906,-61){\Large$\lambda$}
\end{picture}
\caption{$\lambda \ral^\pi \lambda'$, $\pi = (-\frac{3}{2}, \frac{3}{2})$}
\label{lplp}
\end{figure}

Let $\Pi_{k,l}(\lambda)$ be the subset
of pairs $(\alpha, \beta) \in \Pi(\lambda)$ for which $\alpha > -k$ and 
$\beta < l$.
\begin{lemma} Suppose that $\lambda \ral^\pi \lambda'$, and $\lambda \in
\Omega_{k,l}$.  Then $\lambda' \in \Omega_{k,l}$ if and only if 
$\pi \in \Pi_{k,l}(\lambda)$.
\end{lemma}

Our interest in this relation comes because
$\lambda \lra \lambda'$ 
if and only if the corresponding irreducible components 
$\overline{T^*_{X_\lambda}X}$ and $\overline{T^*_{X_{\lambda'}}X}$ of 
the conormal variety $\Lambda$ intersect in
codimension one;  we prove this in \S2 (see Corollary \ref{c1int}).
Lascoux and Schutzenberger studied this relation in 
\cite{LS}; they showed that
$\lambda \ra \lambda'$ if and only if the stalk intersection cohomology group
\[ IH_x^{\lvert \lambda' \rvert - \lvert \lambda \rvert - 1}( 
\overline{X_{\lambda'}})\] is nonzero, where $x$ is 
any point of $X_\lambda$.

\subsection{The quiver category: type A}
In this section, we define the quiver category that describes
perverse sheaves on the Grassmannians $X_{k,l}$.

\begin{defn} A {\em diamond\/} is a 4-tuple 
$(\lambda,\lambda',\lambda'',\lambda''')$ of distinct elements in $\Omega$
satisfying $\lambda \lra \lambda' \lra \lambda'' \lra \lambda''' \lra \lambda$.
\end{defn}

Let $\cA_{k,l}$ be the category whose objects are collections
of $\Q$-vector spaces $V_\lambda$, $\lambda \in \Omega_{k,l}$,
 together with two types of homomorphisms:
\begin{itemize}
\item maps $t_\alpha$ for each $\alpha \in \Hi$, acting on 
$V = \bigoplus V_\lambda$, and preserving this decomposition, and
\item maps $p(\lambda, \lambda')\colon
V_{\lambda'} \to V_\lambda$ 
for every pair $\lambda$, 
$\lambda' \in \Omega_{k,l}$ with $\lambda \lra \lambda'$. 

\end{itemize}  If $\lambda \lra \lambda'$, define
$\mu(\lambda, \lambda') = 1 + p(\lambda, \lambda')p(\lambda', \lambda)$.
These homomorphisms are required to satisfy the following relations:
\begin{enumerate}
\item The maps $t_\alpha$ commute with each other and with the $p$'s.
\item If $(\alpha, \beta)$ is a $\lambda$-pair, then 
$t_\alpha t_\beta|_{V_\lambda} = 1|_{V_\lambda}$.  If $\alpha < -k$ or $\alpha > l$, then 
$t_\alpha = 1$.

\item Suppose that $\lambda \ral^{(\alpha, \beta)} \lambda'$, and 
$(\alpha', \beta')$ is the parent of $(\alpha, \beta)$ in $\Pi(\lambda)$.
Then 
\begin{eqnarray*}\mu(\lambda', \lambda)^{\eta(\beta)}  &=& t^{}_\alpha 
t^{}_{\beta'}|_{V_{\lambda'}},\;\text{and}\\ 
\mu(\lambda, \lambda')^{\eta(\beta)} & = & t^{}_\alpha 
t^{}_{\beta'}|_{V_{\lambda}},
\end{eqnarray*}
where $\eta(\beta) = (-1)^{\beta + \half}$.

\item If $\lambda, \lambda', \lambda'',$ and $\lambda''' \in \Omega_{k,l}$
form a diamond,
then \[p(\lambda'', \lambda')p(\lambda', \lambda) = p(\lambda'', \lambda''')
p(\lambda''', \lambda).\]
If we have a diamond $(\lambda, \lambda', \lambda'', \lambda''')$,
and all elements except $\lambda'''$ are in $\Omega_{k,l}$, then
\[p(\lambda, \lambda')p(\lambda', \lambda'') = 0 = p(\lambda'', \lambda')
p(\lambda',\lambda).\]
\end{enumerate}

Maps between objects of $\cA_{k,l}$ are collections of
maps $V_\lambda \to V'_\lambda$ which are compatible with the 
$t$'s and $p$'s.  Note that the presentation above admits simplifications;
for instance, Lemma \ref{tredund} below shows that the $t_\alpha$ can 
be expressed in terms of the $p$'s.

Objects of $\cA_{k,l}$ are finite dimensional representations of 
the $\Q$-algebra $A_{k,l}$ generated by the $t_\alpha$, 
$p(\lambda, \lambda')$, and commuting idempotents $e_{\lambda'}$ 
representing the projections $\bigoplus V_\lambda \to V_{\lambda'}$, modulo 
the relations above.

\begin{prop} The algebra $A_{k,l}$ is finite dimensional over $\Q$.
\end{prop}
\begin{proof} This follows from Proposition \ref{noval} below.
\end{proof}

\subsection{The main result: type A} 
In \S4 we will define a functor $R\colon \cP_\Lambda(X_{k,l}) \to \cA_{k,l}$
from Schubert-constructible perverse sheaves to quiver representations.
If $R({\mb P}) = \{V_\lambda, t_\alpha, p(\lambda,\lambda')\}$, then 
$V_\lambda$ is the 
vanishing cycles group, or Morse group, of $\mb P$ at a point in the Schubert stratum $X_\lambda$
corresponding to $\lambda$.

\begin{thm} \label{Amain} $R$ is an equivalence of categories.
\end{thm}

The form of the quiver algebra $\cA_{k,l}$ reflects the geometry of
the conormal variety $\Lambda$ associated to the Schubert stratification of $X_{k,l}$. 
The action of $t_\alpha$ represents the monodromy of the vanishing cycles local systems
(which are local systems on the smooth components of $\Lambda$) 
around loops generated by the action of a loop $\gamma_\alpha \in \pi_1(B)$.

The quiver relation (2) describes the kernel of the
homomorphisms $\pi_1(B) \to \pi_1(\wt\Lambda_{X_\lambda})$, where
$\wt\Lambda_{X_\lambda}$ is the smooth part of the conormal variety 
lying over $X_\lambda$.  (3) relates the monodromies of
loops in $\wt\Lambda_{X_\lambda}$ and $\wt\Lambda_{X_{\lambda'}}$
around the intersection of their closures; as we shall see in \S2.5, they intersect
in a divisor if and only if $\lambda \ra \lambda'$.  Finally,
 (4) comes 
from the codimension two strata of $\Lambda$; it is analogous to,
and follows from, a similar relation in the quiver for perverse
sheaves on $\C^2$ stratified by normal crossings.


\subsection{The type D isotropic Grassmannian}
To describe the other space we will consider, let $n = 2k$, and take a 
nondegenerate quadratic form on $\C^n$:
\[Q(x_1, \dots, x_{2k}) = x_1x_{2k} + x_2x_{2k-1} + \dots + x_{k}x_{k+1}.\]
The space $\{V \in X_{k,k} \mid Q|_V = 0\}$ has two connected components;
let $X^s_k$ be the component containing the point $\{x_{k+1},\dots, x_n = 0\}$.  
The group $G^s = SO(n, \C) \subset SL(n, \C)$ of 
transformations preserving $Q$ acts on $X^s_k$, and 
$X^s_k = G^s/(G^s \cap P_{k,k})$.  Throughout this article, we will refer to
this space as the ``type D case", and we will use a superscript ``$s$" to
distinguish structures used in this case.

The group $B^s = G^s \cap B$ is a Borel subgroup of $G^s$.
The Schubert stratification 
of $X^s_k$ is the stratification by $B^s$ orbits, which are just
the spaces $X^s_\lambda = X_\lambda \cap X^s_k$.
Let $\Omega^s \subset \Omega$ be the set of 
partitions whose Young diagrams are symmetric about the diagonal and
 have an even number of squares on the diagonal. We have 
\[\Omega^s =\{\lambda \in \Omega \mid
    \phi_\lambda(\alpha) = -\phi_\lambda(-\alpha)\;\text{for all}\; 
\alpha \in \Hi\;\text{and}\; |\lambda|\; \text{is even}\}.\]
Let $\Omega^s_k = \Omega^s \cap \Omega_{k,k}$;
 then $\lambda \in \Omega^s_k$ if and only if $X^s_\lambda \ne \emptyset$.

\subsection{$\lambda$-pairs: type D}

Again we define a set $\Pi^s(\lambda) \subset \wt{\Pi}(\lambda)$ of 
$\lambda$-pairs for $\lambda \in \Omega^s$. 
\begin{defn}
If $(\alpha, \beta) \in \wt{\Pi}(\lambda)$, 
we let $(\alpha, \beta) \in \Pi^s(\lambda)$ if 
one of the following two conditions holds:
\begin{itemize}
\item $\alpha>0$, $\beta> 0$, and
$(\alpha, \beta) \in \Pi(\lambda)$, or
\item  $\beta>-\alpha>0$, $(\alpha, -\alpha)$ and $(-\beta, \beta)$ are 
both in $\Pi(\lambda)$, $\alpha + \half$ is even, and 
$\sum_{-\alpha<\gamma\le \beta}\phi_\lambda(\gamma)$ = +1.  
\end{itemize}
\end{defn}
(The last equation simply says that $(-\beta, \beta)$ is the parent of
$(\alpha, -\alpha)$ in $\Pi(\lambda)$.) 
Let $\Pi^s_k(\lambda) \subset \Pi_s(\lambda)$ be 
the subset of $(\alpha, \beta)$ for which $\beta < k$.
Give $\Pi^s(\lambda)$ the partial ordering induced from 
$\wt{\Pi}(\lambda)$.

\begin{lemma} For any $\alpha \in \Hi$ there is exactly one type D
$\lambda$-pair in which either $\alpha$ or $-\alpha$ appears.  $\Pi^s(\lambda)$
is a tree with the given partial ordering.
\end{lemma}
\begin{proof} Since the Young diagram of $\lambda$ is symmetric,
we have $(\alpha, \beta) \in \Pi(\lambda)$ if and only if $(-\beta, -\alpha)
\in \Pi(\Lambda)$.  
Elements $(\alpha, \beta)$ where $\alpha$, $\beta$ have the same
sign thus come in pairs; the first part of the definition of $\Pi^s(\lambda)$
just takes the one with positive $\alpha$ and $\beta$.

The remaining elements of $\Pi(\lambda)$ must be of the form 
$(-\alpha, \alpha)$ for $\alpha > 0$.  
These form a sequence of nested intervals
$(\alpha_1,-\alpha_1) < (\alpha_2,-\alpha_2) < \dots.$
Since each of the sets $[0, -\alpha_1)\cap \Hi$ and 
$(-\alpha_i, -\alpha_{i+1}) \cap \Hi$ must be 
a union of pairs in $\Pi(\lambda)$, we see that
$\alpha_i + \half + i$ is even for all $i$.

The second part of the definition 
above thus puts $(\alpha_{2i-1},-\alpha_{2i})$ in $\Pi^s(\lambda)$
for all $i\in \N$. Since any $\alpha$ appears in exactly one 
pair in $\Pi(\lambda)$, the first statement follows. 

The fact that $\Pi^s(\lambda)$ is a tree now follows easily 
from the fact that $\Pi(\lambda)$ is a tree.
\end{proof}

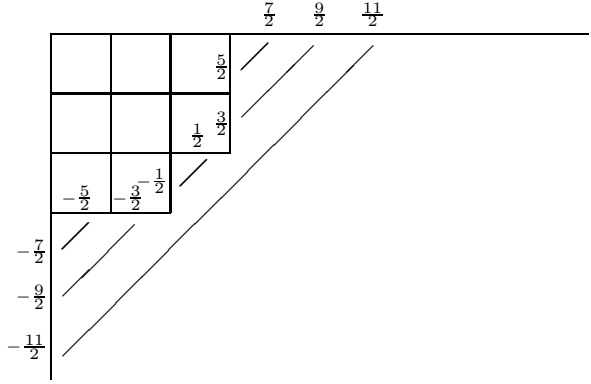
\begin{figure}
\begin{center}

\setlength{\unitlength}{3947sp}%
\begingroup\makeatletter\ifx\SetFigFont\undefined%
\gdef\SetFigFont#1#2#3#4#5{%
  \reset@font\fontsize{#1}{#2pt}%
  \fontfamily{#3}\fontseries{#4}\fontshape{#5}%
  \selectfont}%
\fi\endgroup%
\begin{picture}(3507,2376)(1081,-2023)
\put(1736,-790){$\scriptscriptstyle -\frac{1}{2}$} 
\put(1586,-896){$\scriptscriptstyle -\frac{3}{2}$}
\put(1266,-896){$\scriptscriptstyle -\frac{5}{2}$} 
\put(2071,-501){$\scriptscriptstyle \frac{1}{2}$}
\put(2221,-426){$\scriptscriptstyle \frac{3}{2}$}
\put(2221,-71){$\scriptscriptstyle \frac{5}{2}$} 
\put(2521,259){$\scriptscriptstyle \frac{7}{2}$}
\put(2836,259){$\scriptscriptstyle \frac{9}{2}$}
\put(3141,259){$\scriptscriptstyle \frac{11}{2}$} 
\put(981,-1231){$\scriptscriptstyle -\frac{7}{2}$}
\put(981,-1516){$\scriptscriptstyle -\frac{9}{2}$}
\put(921,-1831){$\scriptscriptstyle -\frac{11}{2}$}
\thinlines
\put(1201,-586){\line( 1, 0){375}}
\put(1201,-211){\line( 1, 0){375}}
\put(1576,-211){\line( 0, 1){375}}
\put(1951,-211){\line( 0, 1){375}}
\put(1951,-586){\line( 0, 1){375}}
\put(1951,-961){\line( 0, 1){375}}
\put(1951,-961){\line( 0, 1){375}}
\put(2326,-586){\line( 0, 1){375}}
\put(1576,-586){\line( 1, 0){375}}
\put(1951,-586){\line( 1, 0){375}}
\put(1576,-961){\line( 1, 0){375}}
\put(1726,-1036){\line(-1,-1){450}}
\put(2401,-361){\line( 1, 1){450}}
\put(1276,-1861){\line( 1, 1){1950}}
\put(1201,-2011){\line( 0, 1){2175}}
\put(1201,164){\line( 1, 0){1125}}
\put(2011,-795){\line( 1, 1){170}}
\put(2396,-60){\line( 1, 1){170}}
\put(1271,-1191){\line( 1, 1){170}}
\put(1201,-961){\line( 1, 0){375}}
\put(1576,-961){\line( 0, 1){750}}
\put(1576,-211){\line( 1, 0){750}}
\put(2326,-211){\line( 0, 1){375}}
\put(2326,164){\line( 1, 0){2250}}
\end{picture}

\end{center}
\caption{Finding $\Pi^s([3,3,2])$.}
\label{Dpairs}
\end{figure}

To give an example, let $\lambda = (3,3,2)$.  Figure \ref{Dpairs} shows that
$(\frac{5}{2}, \frac{7}{2}), (\frac{3}{2}, \frac{9}{2})$ and 
$(-\half, \frac{11}{2})$ are all in $\Pi^s(\lambda)$.

If $(\alpha, \beta) \in \Pi^s(\lambda)$, define a 
second partition $\lambda' \in \Omega$ by switching the signs of 
$\phi_\lambda(\alpha)$, $\phi_\lambda(\beta)$, $\phi_\lambda(-\alpha)$, 
and $\phi_\lambda(-\beta)$.  It is easy to check that $\lambda' \in \Omega^s$.
As before, we write $\lambda \ral^{(\alpha, \beta)} \lambda'$ to 
express this relation.

\subsection{The main result: type D}

Now we define a quiver category $\cA^s_k$ analogous to the 
category $\cA_{k,l}$ from the last section.  Objects are collections
of finite dimensional
vector spaces $V_\lambda$, one for each $\lambda\in \Omega^s_k$,
 together with maps $t_\alpha$ for $\alpha \in \Hi$ and 
$p(\lambda, \lambda')$ for pairs $\lambda \lra \lambda'$, as before.

Once again we define 
$\mu(\lambda, \lambda') = 1 + p(\lambda, \lambda')p(\lambda', \lambda)$
whenever $\lambda \lra \lambda'$.  These maps are then required to satisfy 
the following relations, plus the relation (1) from the definition
of $\cA_{k,l}$:
\begin{itemize}
\item[(2s)] $t_\alpha t_{-\alpha} = 1$ for all $\alpha$, and $t_\alpha = 1$ if 
$\alpha> k$.  
If $(\alpha, \beta)$ is a $\lambda$-pair, then
$t_\alpha t_\beta|_{V_\lambda} = 1|_{V_\lambda}$.
\item[(3s)] Suppose that $\lambda \ral^{(\alpha, \beta)} \lambda'$,
 and $(\alpha', \beta') \in \Pi^s(\lambda)$ is the parent of $(\alpha, \beta)$.
If $(\alpha, \beta) < (-\beta', -\alpha')$, then let $\zeta = -\alpha'$;
otherwise let $\zeta = \beta'$.
Then
\begin{eqnarray*}\mu(\lambda', \lambda)^{\eta(\beta)} &=& t^{}_\alpha 
t^{}_{\zeta}|_{V_{\lambda'}},\;\text{and}\\ 
\mu(\lambda, \lambda')^{\eta(\beta)} & = & t^{}_\alpha 
t^{}_{\zeta}|_{V_{\lambda}},
\end{eqnarray*}
where $\eta(\beta) = (-1)^{\beta + \half}$.
(Note that because of relation 2s, the second equation can
also be written
$\mu(\lambda, \lambda')^{\eta(\beta)}  =  t^{}_\alpha 
t^{}_{\beta'}|_{V_{\lambda}}$.)
\item[(4s)] The first sentence of relation (4) from $\cA_{k,l}$ holds, 
but the second part is modified as follows:
If we have a triple $\lambda \lra \lambda' \lra \lambda''$ in $\Omega^s_k$ 
for which either (a) there is 
a diamond $(\lambda, \lambda', \lambda'', \lambda''')$,
where $\lambda''' \in \Omega^s \setminus \Omega^s_k$ 
{\bf or} (b) $\lambda \ral^\pi \lambda' \ral^{\pi'}
\lambda''$ with $\pi = (\alpha, \beta)$ a $\lambda$-pair with $\alpha<0$ and 
$\pi' \notin \Pi^s(\lambda)$, then 
\[p(\lambda, \lambda')p(\lambda', \lambda'') = 0 = p(\lambda'', \lambda')
p(\lambda',\lambda).\]
\end{itemize}
Here the definition of a diamond is the same as in \S1.3: a $4$-tuple
$(\lambda,\lambda',\lambda'',\lambda''')$ of distinct elements in $\Omega^s$
satisfying $\lambda \lra \lambda' \lra \lambda'' \lra \lambda''' \lra \lambda$.

Define an algebra $A^s_k$ in the same way as before, so that $\cA_k^s$
is the category of finite dimensional representations of $A^s_k$. 

\begin{prop} The algebra $A^s_k$ is finite dimensional over $\Q$.
\end{prop}
\begin{proof} This follows from Proposition \ref{noval} below.
\end{proof}

We will define a functor $R\colon \cP_\Lambda(X^s_k) \to \cA^s_k$
just as in the type $A$ case.  Our main result for the type $D$ 
Grassmannian is:

\begin{thm} \label{Dmain} $R$ is an equivalence of categories.
\end{thm}

The definition of $R$ and the proof of Theorem \ref{Dmain} will
be given in \S4.

\subsection{Simple representations and finite dimensionality}
In this section we prove that the algebras $A = A_{k,l}$ and 
$A = A^s_k$ are
finite dimensional, and describe their irreducible representations.

\begin{lemma} \label{tredund} Take 
 $\gamma \in \Hi$, and assume that $\gamma>0$ in the type D case. Then
we have \[(t_\gamma|_{V_\lambda})^{\phi_\lambda(\gamma)} = 
\prod \mu(\lambda, \lambda')^{\eta(\beta)},\]
where the product is over all $\lambda'\in \Omega_{k,l}$ 
(resp. $\Omega^s_k$) for which 
$\lambda \ral^{(\alpha, \beta)} \lambda'$ for some $\lambda$-pair with
$\alpha \le \gamma \le \beta$.
\end{lemma}
\begin{proof}
  Using the quiver relations 2 and 3 (or 2s and 3s), we see that
if $\lambda \ral^{(\alpha, \beta)} \lambda'$, and 
$(\alpha', \beta')$ is the parent of $(\alpha, \beta)$ we have
$\mu(\lambda, \lambda')^{\eta(\beta)} 
= t^{}_\alpha t^{-1}_{\alpha'}$ on $V_\lambda$.
The product above is thus a telescoping product, since $t^{}_\alpha = 1$
for all sufficiently negative $\alpha$.  The telescope starts with 
$t_\gamma^{\pm 1}$, where the sign is determined by $\phi_\lambda(\gamma)$,
i.e.\ by whether 
$\gamma$ appears first or second in a $\lambda$-pair.
\end{proof}

\begin{prop} \label{flippi1} If $\lambda' \ra \lambda$ then
$\mu(\lambda, \lambda')$ and $\mu(\lambda', \lambda)$ are both 
unipotent in $A$.  
\end{prop}
\begin{proof} First note that $\mu(\lambda,\lambda')$ is unipotent if 
and only if $\mu(\lambda',\lambda)$ is, since they are of the
form $1+pq$ and $1+qp$, respectively.

We proceed by downward induction in $\lambda$.  If $\lambda$
is the maximal element in $\Omega_{k,l}$ (or $\Omega^s_k$), then
the quiver relation 3 implies
$\mu(\lambda, \lambda') = 1$, since all $t_\alpha$ act as 
the identity on $V_\lambda$.

Now suppose the proposition holds for all $\lambda > \lambda_0$,
and take some $\lambda' \ra \lambda_0$.
Lemma \ref{tredund} and the quiver condition 3 can be used to express 
$\mu(\lambda_0, \lambda')$ as a product of powers (positive and negative) of $\mu(\lambda_0, \lambda'')$
for $\lambda_0 \ra \lambda''$; these are unipotent by the inductive assumption, and 
they commute by quiver relation 1.  The unipotence of 
$\mu(\lambda_0, \lambda')$ and $\mu(\lambda', \lambda_0)$ now follow.
\end{proof}

\begin{lemma} \label{diam}
Suppose that $\lambda \ra \lambda_1$, $\lambda \ra \lambda_2$
and $\lambda_1\ne \lambda_2$
 (in either $\Omega$ or $\Omega^s$).  Then we can complete this
arrangement to a diamond, i.e.\
there exists a
 $\lambda'\in \Omega$ (resp. $\Omega^s$) with $\lambda' \ne \lambda$
for which $\lambda_1 \lra \lambda' \lra \lambda_2$.
Furthermore, for any such $\lambda'$, 
either $\lambda'>\lambda_1$ or $\lambda'>\lambda_2$.
\end{lemma}
\begin{proof}
Consider the case of $\Omega$; the argument for $\Omega^s$ is similar and
we will omit it.  Suppose that 
$\lambda \ral^{\pi_1} \lambda_1$ and $\lambda \ral^{\pi_2} \lambda_2$
Recall the tree structure on the set $\Pi(\lambda)$ of $\lambda$-pairs.
If neither $\pi_1$ or $\pi_2$ is the parent of the other, then we
have $\pi_1 \in \Pi(\lambda_2)$ and $\pi_2 \in \Pi(\lambda_1)$, 
and there is a $\lambda'$ with $\lambda_1 \ral^{\pi_2} \lambda'$ and
$\lambda_2\ral^{\pi_1}\lambda'$.

On the other hand, if $\pi_2 = (\alpha_2, \beta_2)$ is the parent of 
$\pi_1 = (\alpha_1, \beta_1)$, then $(\alpha_2, \alpha_1)$ and 
$(\beta_1, \beta_2)$ are in $\Pi(\lambda_1)$, and 
we get $\lambda_l, \lambda_r$ for which
 $\lambda_1 \ral^{(\alpha_2, \alpha_1)} \lambda_l \ral^{(\beta_1, \beta_2)} 
\lambda_2$ and
$\lambda_1 \ral^{(\beta_1, \beta_2)} \lambda_r \ral^{(\alpha_2, \alpha_1)} 
\lambda_2$.

For the last statement, just check that
in both cases above we have found the only possible $\lambda'\ne \lambda$ for which
$\lambda_1 \lra \lambda' \lra \lambda_2$.
\end{proof}

Since $p(\lambda, \lambda')p(\lambda'', \lambda''') =0$ unless 
$\lambda' = \lambda''$, a nonzero monomial in the $p$'s can be 
seen as a path in the graph $\Gamma$ whose nodes are elements of
$\Omega_{k,l}$ (or $\Omega^s_k$), 
with edges are given by the relation ``$\lra$". 
Let the idempotent $e_\lambda$ be
the monomial corresponding to the trivial path
with only one node $\lambda$.  
Let $(\lambda_1, \lambda_2, \dots, \lambda_j)$ be a path in $\Gamma$.
We will say that $\lambda_i$ is a {\em valley\/} if  
$\lambda_{i-1} > \lambda_i$ and $\lambda_{i+1} > \lambda_i$. 

\begin{prop} \label{noval}
The algebra $A$ is spanned as a vector space
over $\C$ by the monomials without valleys.  A monomial corresponding to 
a path of length $l$ is a $\Z$-linear combination of monomials without 
valleys, all coming from paths of length $\ge l$.
\end{prop}

\begin{proof}
First, note that as a consequence of Lemma \ref{tredund}, we only need
to show that monomials in the $p$'s can be expressed in terms of
monomials without valleys.

Consider the monomial $m$ with path 
$(\lambda_1, \lambda_2, \dots , \lambda_j)$, and let $\lambda_i$ be
a valley.  If $\lambda_{i-1} \ne \lambda_{i+1}$, then
we can apply Lemma \ref{diam} to obtain another expression for 
$m$ as a monomial of the same length without a valley in the $i$th place.  

If $\lambda_{i-1} = \lambda _{i+1}$, then we have
$p(\lambda_{i-1}, \lambda_i)p(\lambda_i, \lambda_{i-1})=
\mu(\lambda_{i-1}, \lambda_i) - 1$.
Now apply Lemma \ref{tredund} to express 
$\mu(\lambda_{i-1}, \lambda_i)$ as a product of terms 
$\mu(\lambda_{i-1},\lambda)^{\pm 1}$ for 
$\lambda > \lambda_{i-1}$.  Since all these terms are unipotent,
we can use the substitution $(1 + pq)^{-1} = \sum (-pq)^k$ to get an
expression for $\mu(\lambda_{i-1}, \lambda_i)$ 
as a linear combination of monomials with length $\ge 2$ 
and which only visit nodes $\lambda \ge \lambda_{i-1}$.

In both cases, the only new valleys created are above $\lambda_i$.  
So we can apply this process repeatedly, first getting rid of
all valleys for which $|\lambda_i| = 0$, then for $|\lambda_i| =1$,
and so on.  This process will terminate, since any monomials
whose paths go outside $\Omega_{k,l}$ (resp. $\Omega^s_k$) are zero.
\end{proof}

The finite dimensionality of $A_{k,l}$ and $A^s_k$ follows:
any path without valleys is a composition of an increasing
path with a decreasing path, and there are clearly only 
finitely many of these.  As a further consequence, we obtain
the following description of the irreducible representations of the algebras
$A_{k,l}$ and $A^s_k$, which we will need in the final step 
of our proof of theorems $\ref{Amain}$ and $\ref{Dmain}$.

\begin{thm} \label{simples}
All irreducible representations of $A_{k,l}$ and $A^s_k$
are one dimensional.  They are in one-to-one correspondence
with elements $\Omega_{k,l}$ (resp. $\Omega^s_k$).
\end{thm}
\begin{proof}
Suppose $V = \{V_\lambda, t_\alpha, p(\lambda,\lambda')\}$ 
is an irreducible representation.  Choose a nonzero
vector $v \in V_\lambda$.  Let $(\lambda_1, \dots, \lambda_j)$ be
the longest path for which the corresponding monomial $m$ acts
nontrivially on $v$; one exists because of Proposition \ref{noval}.  
Then for any $\lambda \lra \lambda_j$,
we have $p(\lambda, \lambda_j)m\cdot v = 0$.  If $V' = 
\{V'_\lambda, t'_\alpha,p'(\lambda,\lambda')\}$
is the irreducible representation for which $V'_{\lambda_j}$ is one
dimensional and all other $V'_\lambda$ vanish, then there is a map $V'\to V$
given by sending a generator of $V'_{\lambda_j}$ to $m\cdot v$.  
This contradicts the irreducibility of $V$ unless the path was trivial
and $V' = V$.
\end{proof}

\begin{rem} 
Combining Theorems \ref{simples}, \ref{Amain}, and \ref{Dmain},
we see that simple perverse sheaves in $\cP_\Lambda(X_{k,l})$ and
$\cP_\Lambda(X^s_k)$ have nonzero vanishing cycle groups at only
one stratum.  This was proved in the type A case by Bressler, Finkelberg
and Lunts \cite{BFL}, and in the type D case by Boe and Fu \cite{BF}.
\end{rem}

\section{Conormal geometry of the Grassmannian}
\label{mlgeom}
In this section we study the geometry of the conormal variety $\Lambda$
to the Schubert stratifications of $X = X_{k,l}$ and $X = X^s_k$.  The
key facts, in both cases, are:
\begin{itemize}
\item The Borel group $B$ acts on $\Lambda$ with finitely many orbits.
\item The stabilizers of the $B$-action on $\Lambda$ are connected, 
so the fundamental groups of orbits are quotients of $\pi_1(B)$.
\item Near a point in a codimension one or two orbit $\Lambda$ has only
normal crossings singularities.
\end{itemize}

More precisely we focus on 
the fiber $M_\lambda$ of $\Lambda$ over a torus-fixed point $W_\lambda$ in 
$X_\lambda$.  The stabilizer $B_\lambda = B_{W_\lambda}$ acts on this fiber 
with finitely many orbits (Proposition \ref{Mlorbs}).  The $B$-orbits
of $\Lambda$ lying over $X_\lambda$ are isomorphic to $X_\lambda \times O$,
for $O\subset M_\lambda$ a $B_\lambda$-orbit.

We begin with the case $X = X_{k,l}$, and finish with an outline
of the differences in the type D case.

\subsection{Normal and conormal coordinates} \label{ncc} We first introduce
the coordinate systems we will use to describe the geometry of
$X$ and $\Lambda \subset T^*X$.  

  Given a partition $\lambda \in 
\Omega_{k,l}$, let $W_\lambda$ be the unique point of $X_\lambda$ 
which is fixed by the torus $T \subset B$ of diagonal matrices.
More explicitly, $W_\lambda$ is spanned by $\{\,\be_i \mid i\in I\,\},$ where
\begin{eqnarray*}
I & =  & \{\,\lambda_k +1,\, \lambda_{k-1} + 2,\, \dots, 
\lambda_1 + k\,\}\\
& = & \{\, \alpha + k + {\textstyle\half} \mid \phi_\lambda(\alpha) = -1 \;\text{and}\;
-k < \alpha < l\,\}.
\end{eqnarray*}
Let $I' = \{1, \dots, n\} \setminus I$.  For the rest of this section the
partition $\lambda$ will be fixed and we will put $W = W_\lambda$.

Let $G = SL(n, \C)$, and let $P = G_W$ be the stabilizer of $W$, with Lie
algebra $\mathfrak{p}$.
The infinitesimal action of $\mathfrak g$ on $X$ induces an isomorphism
$T_WX \cong \mathfrak{g/p}$.  Since 
$\mathfrak{p} = \{g\in \mathfrak{g}\mid g(W) \subset W\}$,
we get identifications $T^{}_WX\cong \Hom(W, \C^n/W)$ and
$T^*_WX \cong \Hom(\C^n/W, W)$.  

We represent elements of $\Hom(\C^n/W, W)$ by $k\times l$ matrices,
where we take $\{\be_i\}_{i\in I}$,
in order of {\em decreasing\/} $i$, as a basis for $W$, and the residues of
$\{\be_i\}_{i \in I'}$, taken in increasing order, as 
a basis $\C^n/W$.  
Let $R$ be the rectangle $\{1,\dots,k\}\times\{1,\dots,l\}$. 
For $A\in \Hom(\C^n/W, W)$,
define the {\em support\/} $\Supp A \subset R$ of $A$ to be the 
set of $(i, j) \in R$ such that $A_{ij} \ne 0$.  For a subset 
$S \subset R$, let $\C^S = \{\,A \mid \Supp A \subset S\,\}$.  We use the
usual convention for matrix coordinates so $(i, j)$ 
denotes a point in the $i$th row from the top and the $j$th column from the
left.  We will abuse notation slightly and represent elements 
of $\Hom(W, \C^n/W)$ by $k \times l$ matrices also, using the 
duality given by the basis of elementary matrices.

Let $M_\lambda = T^*_WX \cap T^*_{X_\lambda}X$, the fiber of $\Lambda$
over $W$. Put $\Rb = R\setminus \Delta(\lambda)$, where 
$\Delta(\lambda)$ is the Young diagram of $\lambda$, as introduced in 
\S1.1.

\begin{lemma} \label{abcde}
With the above conventions, $M_\lambda = \C^\Rb.$
\end{lemma}

Before proving this, we need to 
define a map relating the $(i,j)$ matrix coordinates with coordinates
in $\C^n$.  Let 
$w^{}_\lambda\colon \wt{\Pi}(\lambda)\to \N \times \N$ be given by
$w^{}_\lambda(\alpha, \beta) = (w_1(\alpha),w_2(\beta))$,
where \[w_1(\alpha) = \#\{\, \alpha'\in \Hi \mid 
\alpha'\ge\alpha\;\,\text{and}\;\,
\phi_\lambda(\alpha') = -1\,\},\,\text{and}\] 
 \[w_2(\beta) = \# \{\,\beta' \in \Hi \mid \beta'\le \beta
\;\, \text{and}\;\, \phi_\lambda(\beta') = +1\,\}.\] 
Most of the time, the element $\lambda$ will be fixed and we will
drop the subscript $\lambda$.
In terms of the Young diagram pictures,  
$w(\alpha, \beta)$ is the point in the same row as 
the vertical segment on the boundary of 
$\Delta(\lambda)$ indexed by $\alpha$ and the same column as
 the horizontal segment indexed by $\beta$.

\begin {proof}[Proof of Lemma \ref{abcde}]
 It is enough to show that $T_WX_\lambda= \C^{\Delta(\lambda)}$.
 Clearly $T_WX_\lambda$ is the image of the map
 $\mathfrak{b} \to \mathfrak{g/p} \cong \Hom(W, \C^n/W)$.
Diagonal elements of $\mathfrak{b}$
are in the kernel, and the elementary matrix $E_{ij}$
$(i \ne j)$ is in $\mathfrak{p}$ unless $i \in I'$, $j\in I$, 
in which case it maps to the elementary matrix 
with coordinate $w(j -k - \half, i-k-\half)$ in
 $\Hom(W, \C^n/W)$. The result now follows from the following lemma.
\end{proof}

\begin{lemma}
\label{xxx}
If $(\alpha, \beta) \in \wt{\Pi}(\lambda)$,
we have $w(\alpha, \beta) \in \Delta(\lambda)$
if and only if $\alpha > \beta$.
\end{lemma}

Next define a map $\tilde\epsilon = \tilde\epsilon_\lambda \colon \Hom(W, \C^n/W) \to X$
by identifying 
$\C^n$ with $W \oplus \C^n/W$ using the basis $\{\be_i\}$ and 
letting $\tilde\epsilon(A)$ be the graph of $A$.  This embeds 
$T_WX$ as a tubular neighborhood of $X_\lambda$.

Let $\bar{M}_\lambda = \C^{\Delta(\lambda)}$, so we have a splitting
$\Hom(W,\C^n/W) = M_\lambda \oplus \bar{M}_\lambda$.  This gives an 
inclusion $M^*_\lambda \subset \Hom(W, \C^n/W)^* = \Hom(\C^n/W,W)$.
Let $\epsilon\colon M^*_\lambda \to X$ be the restriction of 
$\tilde\epsilon$.  The next result shows that $\epsilon$ is the 
inclusion of a normal slice to $X_\lambda$, and the stratification 
in $\tilde\epsilon(T_WX)$ is the product of $X_\lambda$ with the 
stratification in the normal slice.

\begin{prop} \label{eptilde} If $q\colon \Hom(W, \C^n/W) \to M^*_\lambda$ 
is the projection
map, then $\tilde\epsilon^{-1}(X_{\lambda'}) = q^{-1}\epsilon^{-1}
(X_{\lambda'})$
for any $\lambda'\ge \lambda$.
\end{prop}
\begin{proof} Suppose $A \in \Hom(W, \C^n/W)$ has $A_{xy} = c \ne 0$ for some
$(x,y) \in \Delta(\lambda)$.  If $(i-k-\half,j-k-\half) = 
w^{-1}(x,y)$, then
$\be_i$, $\be_j$ are the standard basis elements corresponding to the
column and row of the square $(x,y)$.  By Lemma \ref{xxx}, we have 
$i>j$.  Then acting on $\tilde\epsilon(A)$ 
by the matrix $I-cE_{ji} \in B$ kills the entry at $(x,y)$.
Repeating this argument shows that $\tilde\epsilon(A)$ and 
$\epsilon\circ q(A)$ always lie in the same Schubert cell.
\end{proof}

\subsection{Actions and orbits}
Keeping the notation $W = W_\lambda$, 
let $B_\lambda = B_{W} \subset B$ be the stabilizer of
$W$; it acts on $T^*X$ fixing $T^*_WX$ and  
$\Lambda$, so it acts on $M_\lambda$.  To describe
this action, take an element $g\in B_\lambda$.  It induces 
endomorphisms $g_1$ and $g_2$ on $W$ and on $\C/W$. 
If $A \in M_\lambda$, then
$g \cdot A = g_1 A g_2^{-1}$.  In other words,
the action is generated by the elementary row
and column operations, where row $i$ can only be added to row $i'$ for $i < i'$
and column $j$ to column $j'$ for $j < j'$.  

The action of $B_\lambda$ on the dual space $M_\lambda^*$ can 
be expressed similarly, in terms of ``truncated row and column
operations":  a row can be added to a row above and a column 
to a column to the left, but anything appearing inside $\Delta(\lambda)$ 
must be discarded.
Note that both these actions include multiplication by scalars,
so all the orbits are conical.  

\begin{prop} \label{Lambda orbits}
There is a one-to-one correspondence sending $B$-orbits
$\cO \subset \Lambda$ to pairs $(\lambda, O)$ where 
$\lambda \in \Omega_{k,l}$ (or $\Omega^s_k$) and $O$ is a 
$B_\lambda$ orbit in $M_\lambda$.  It is given by choosing
$\lambda$ so that $\cO$ lies over the Schubert cell $X_\lambda$,
and letting $O = M_\lambda \cap \cO$.
\end{prop}

\newcommand{\zo}{{$0$-$1$} }
We will call a matrix $A$ a ``\zo matrix" if all its entries $A_{ij}$
are 0 or 1, and
each row and column has at most one nonzero entry.
The following proposition follows easily from the corresponding
result for square matrices.

\begin{prop} \label{Mlorbs} $B_\lambda$ acts with finitely many orbits on both
$M_\lambda$ and $M^*_\lambda$.  In both cases, any orbit contains
a unique point given by a \zo matrix.  The orbit a matrix
$A$ belongs to determines and 
is determined by the ranks of the submatrices
$A_\rho$, where $\rho \subset R$ runs over all rectangles
$\rho \subset R_{k,l}$ with $(1,1)\in \rho$ (for $M_\lambda$) and
rectangles $\rho \subset \Rb$ with $(k,l)\in \rho$ (for $M^*_\lambda$).
If $A$ and $A'$ are \zo matrices, $B_\lambda A' \subset 
\overline{B_\lambda A} \iff
\rank A'_\rho \le \rank A_\rho$ for all $\rho$.
\end{prop}
Note that $\rank A_\rho = \#(\rho \cap \Supp(A))$ if $A$ is a \zo matrix.

\begin{cor} The orbit stratifications of $M_\lambda$ and $M^*_\lambda$
are dual, i.e.\ there is a bijective correspondence $O \to O^*$
from the orbits of $M_\lambda$ to the orbits of $M_\lambda^*$ so that
the closure $\overline{O}$ is the dual cone to $\overline{O^*}$.
\end{cor}
\begin{proof} The dual cone to $\overline{O}$ is a $B_\lambda$-invariant
irreducible variety, and hence the closure of an orbit $O^*$.
\end{proof}

The maps $\tilde\epsilon$ and $\epsilon$ 
are not $B_\lambda$-equivariant, but they are equivariant under the maximal torus in 
$B_\lambda$.  We have the following result:

\begin{prop} \label{hjk} For every $\lambda$ the set $\epsilon^{-1}(X_\lambda)$
is a union of $B_\lambda$-orbits.
\end{prop}
\begin{proof} In \cite{BF}, Boe and Fu show that the 
$\epsilon^{-1}(\overline{X_\lambda})$ are cut out by conditions on the ranks
of submatrices containing $(k,l)$ and touching, but not crossing, the boundary
of $\Delta(\lambda)$.  Since the orbits are determined by the ranks of all
rectangular submatrices containing $(k,l)$, the orbit decomposition is finer.
\end{proof}

\begin{figure}
\begin{center}
\leavevmode
\hbox{
\epsfxsize=4.5in
\epsffile{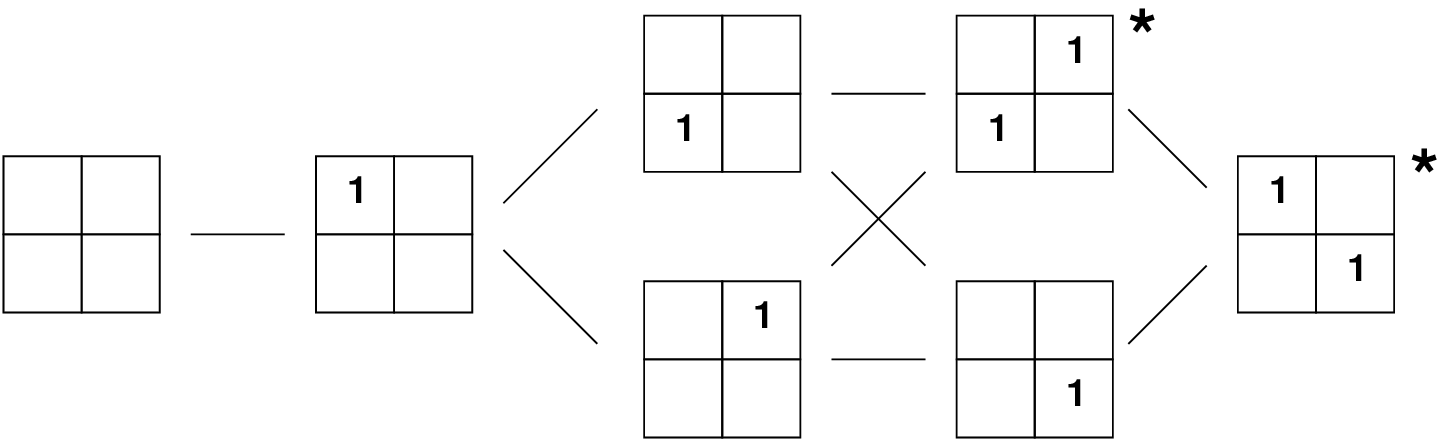}}
\end{center}
\caption{The orbit structure of $M^*_\lambda$, $k = l = 2$, $\lambda = 
\emptyset$.}
\label{orbfig}
\end{figure}

In general the two decompositions of $M^*_\lambda$ are not the same.
For instance, take $k=l=2$, $\lambda=\emptyset$
 the zero partition.  Figure \ref{orbfig} shows the \zo matrices
of the orbits of $M^*_\lambda$ (omitting the zeros).  The lines 
give the codimension one closure relations, with the larger orbit
placed to the right of the smaller one.  Each $\epsilon^{-1}(X_\lambda)$
consists of a single $B_\lambda$ orbit, except for $\lambda = (2, 2)$, 
when it is the union of the two orbits labeled with asterisks.  

The following lemma characterizes which $B_\lambda$ orbits do correspond to Schubert cells.

\begin{lemma} \label{Mstorbs} Take a \zo matrix $A$, considered as an
element of $M^*_\lambda$.  Then $B_\lambda A$ is the 
(unique) open orbit in some $\epsilon^{-1}(X_{\lambda'})$ if and only if
for each $(i, j)$, $(r, s) \in \Supp A$ with $i<r$, $j>s$ we have $(i, s) \in \Delta(\lambda)$.
\end{lemma}
\begin{proof}  It is not hard to see that $\overline{O_A}$
can be defined by using only 
conditions on the ranks of submatrices $A_\rho$ where
$\rho \subset \Rb$ is a rectangle with corners at $(k,l)$
and $(i+1, s+1)$ and 
where $(i,j), (r, s) \in \Supp A$ satisfy $i<r$ and $j>s$.  
The assumption of this lemma assures that such rectangles touch the 
boundary of the Young diagram $\Delta(\lambda)$, so the 
result follows from the characterization 
of normal slices to Schubert cells in \cite{BF}.  
\end{proof}

One case in particular will be important.  Say $E = E_{ij}$ is the 
elementary matrix supported at $(i,j) \in \Rb$, and put
$(\alpha, \beta) = w^{-1}(i,j)$.  The previous lemma shows that
$B_\lambda E \subset M^*_\lambda$ is the open orbit in 
$\epsilon^{-1}(X_{\lambda'})$ for some $\lambda'$.

\begin{lemma} \label{elemmat} We have 
$\phi_{\lambda'}(\alpha) = +1$,
$\phi_{\lambda'}(\beta) = -1$, and 
$\phi_\lambda(\gamma) = \phi_{\lambda'}(\gamma)$ for all 
$\gamma \ne \alpha, \beta$.
\end{lemma}

\begin{proof} $\epsilon(0) = W_\lambda$ is spanned by the vectors
\[\{\be_m \mid 1 \le m \le n,\; \text{and}\; 
\phi_\lambda(m - k - {\textstyle\half}) = -1\} .\]
A basis for $\epsilon(E)$ can be obtained by replacing 
$\be_{\alpha + k + \half}$ by 
$\be_{\alpha + k +\half} + \be_{\beta + k + \half}$.
Since $\alpha<\beta$ (Lemma \ref{xxx}), there is a $g \in B$ for which 
$g\cdot \epsilon(E) = W_{\lambda'}$ is spanned by the same basis as \
$W_\lambda$ with 
$\be_{\alpha+k+\half}$ replaced by $\be_{\beta+k+\half}$.
\end{proof}

\begin{rem} We will see in \S\ref{conesec} that if 
$(\alpha,\beta)$ is a $\lambda$-pair, the dual orbit
$(B_\lambda E)^*$ has codimension one in $M_\lambda$; it is the fiber over
$W_\lambda$ of a codimension one orbit of $\Lambda$ where 
$\overline{T^*_{X_\lambda}X}$ and $\overline{T^*_{X_{\lambda'}}X}$
intersect.
\end{rem}

\subsection{Orbit structure of $M_\lambda$}\label{orbstrs}
In this section we give some general results on the geometry of the
$B_\lambda$-orbits on $M_\lambda$, and the closure relations
between them.  In the end we will only need to understand
the orbits with codimension $\le 2$, but it will be convenient to 
study the problem in general first.

Put a partial order on $\N \times \N$ by saying $(i,j) \le (i',j')$ whenever 
$i\le i'$ and $j\le j'$.  Thus $w$ is an order-preserving map, with the 
order on $\wt\Pi(\lambda)$ from \S1.2.

\begin{figure}
\begin{center}
\leavevmode
\hbox{
\epsfxsize=1.6in
\epsffile{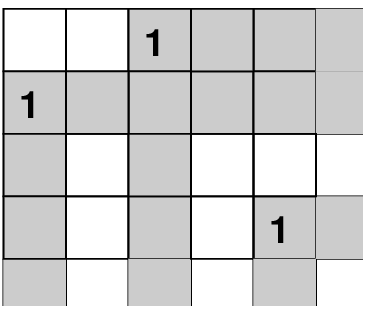}}
\end{center}
\caption{The set $\tau(A)$.}
\label{tauA}
\end{figure}

For a \zo matrix $A \in M_\lambda$, denote the orbit 
$B_\lambda A$ by $O_A$.
Define $\tau(A) \subset \Rb$ to be the set of points $(i,j)$
 for which there is a point $(i', j') \in \Supp A$ with $(i', j') \le (i,j)$
and either $i = i'$ or $j = j'$; i.e., $\tau(A)$ consists of all squares
for which there is a $1$ either above or to the left; see figure \ref{tauA}.

\begin{prop} \label{orbdim}
Suppose $A\in M_\lambda$
is a \zo matrix.  Then we have 
\[T_A(O_A) = \C^{\tau(A)}, 
\]
where we identify the tangent space to the conical orbit $O_A = B_\lambda\cdot A$
as a subspace of $M_\lambda$. 
In particular, we have
$\dim O_A = \#\tau(A)$. \end{prop}

It will be useful to reconstruct the matrix $A$ from $\tau(A)$:
\begin{prop}
\label{tatoa}  A \zo matrix $A$ can be recovered from $\tau(A)$ by 
the following procedure:  begin with $A = 0$, and $S = \tau(A)$.  If
$S=\emptyset$ then halt; otherwise
choose a minimal element $(i,j)$ of $S$ and put $A_{ij}=1$.  Discard all 
points of $S$ in the same row or column as $(i,j)$ and repeat.
\end{prop} 

\begin{figure}
\begin{center}
\leavevmode
\hbox{
\epsfxsize=2in
\epsffile{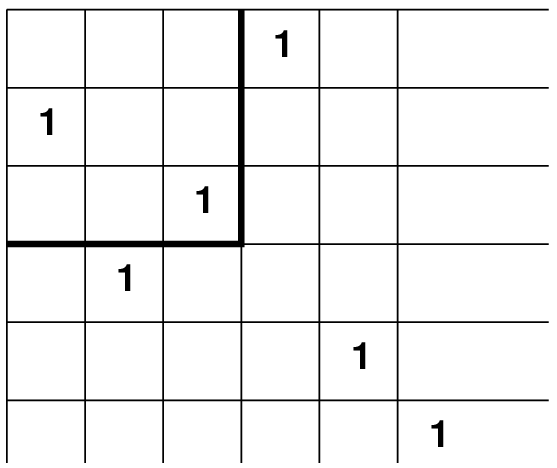}}
\end{center}
\caption{Completing $A$ to $\hat{A}$.}
\label{Ahatfig}
\end{figure}

Given a \zo matrix $A$, we can complete $A$ to an ``infinite \zo matrix"
$\hat{A}$ with entries for all  $(i,j) \in \N^2 \setminus \Delta(\lambda)$
by applying this algorithm to the set 
\[\tau(A) \cup \{(i, j) \mid i>k \;\text{or}\; j > l\}.\]
It follows that $\hat{A}_{ij} = A_{ij}$ for all $(i,j) \in \Rb$.  
Figure \ref{Ahatfig} shows an example where $k = l = 3$, 
$\lambda =\emptyset$.

Since $\hat{A}$ has exactly
one $1$ in each row and column, it can be thought of as a permutation 
$\sigma_A\colon \N \to \N$: let $\sigma_A(i) = j$ whenever $\hat A_{ij} = 1$.  
Then $\sigma$ defines an injective 
map from the set of \zo matrices on $\Rb$ to the group 
$S_\infty = \bigcup S_n$
of permutations of $\N$ which are eventually the identity.

\begin{lemma} \label{codislen} We have 
\[\mathop{\rm len}(\sigma) = kl - \dim(O_A),\] where $\mathop{\rm len}$
is the length function on $S_\infty$ considered as a limit of Coxeter 
groups.  Furthermore, we have $O_{A'} \subset \overline{O_A}$ for 
\zo matrices $A, A'$ if and only if the corresponding 
permutations satisfy $\sigma \le \sigma'$ in the
Bruhat order on $S_\infty$.
\end{lemma}
\begin{proof} The first statement follows since $(i,j) \notin \tau(A)
\iff i < \sigma^{-1}(j)$ and $\sigma(i) > j$, whereas 
$\mathop{\rm len}(\sigma)$ is the number of pairs $(i, j')$ with $i< j'$ and $\sigma(i)>\sigma(j')$.  

For the second statement, we use the following 
description of the Bruhat order on $S_\infty$ 
(see \cite{F}, page 173).  Given $\sigma \in S_\infty$ and $p, q \in \N$, 
\[r_\sigma(p, q) = \#\{i \le p \mid \sigma(i)\le q\}.\] 
Then $\sigma \le \sigma'$ if and only if 
$r_\sigma(p, q) \ge r_{\sigma'}(p,q)$ for all $p, q$.
If $\sigma$ is associated to the matrix $\hat{A}$, then 
$r_\sigma(p, q)$ is just
the rank of the submatrix $\hat{A}(p,q)$
 with corners at $(1,1)$ and $(p, q)$, so the 
``if" part of the lemma is clear.  

For the other direction, note that for points $(p, q) \notin \Rb$ 
we can recursively find 
$r^{}_A(a, b) = \rank \hat{A}(p, q)$: set $r_A(a, 0) = r_A(0, b) = 0$ for all 
$a,b$.  Then
\[r^{}_A(p, q) = \min(r^{}_A(p-1, q), r^{}_A(p, q-1)) + 1.\]
Thus knowing $r_A(p,q) \le r^{}_{A'}(p,q)$ holds for $(p,q) \in R_{k,l}$ 
will imply the same inequality for
all $(p,q)$, giving the ``only if" part.
\end{proof}

\begin{figure}
\begin{center}
\leavevmode
\hbox{
\epsfxsize=3.6in
\epsffile{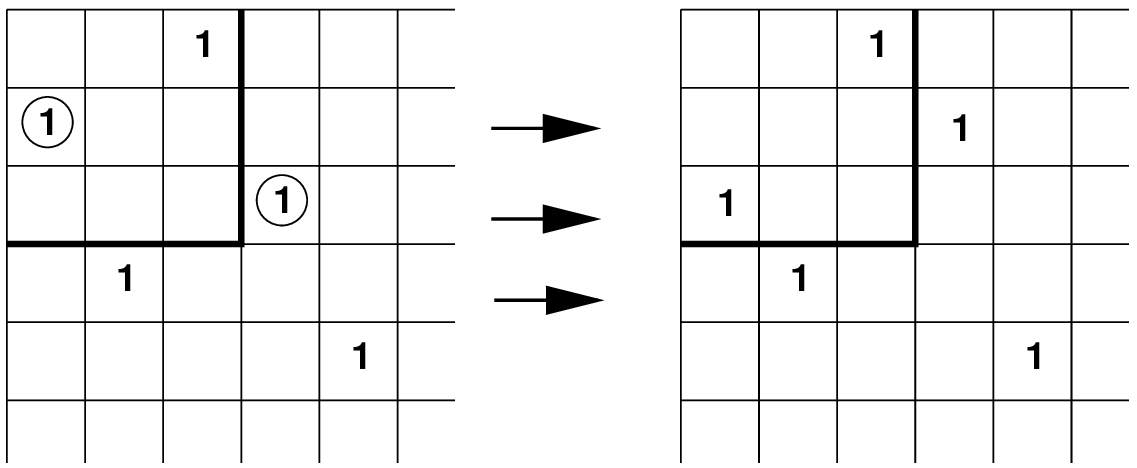}}
\end{center}
\caption{}
\label{fundmove}
\end{figure}

Next we describe a procedure which given an orbit, gives
a codimension one orbit contained in its closure.  Let $A \in M_\lambda$
be a \zo matrix, and take a point $(i,j)$ in its support.
Take a minimal point $(r, s) > (i,j)$ with $\hat{A}_{rs} = 1$, and  
define $A'$ by
\begin{itemize}
\item $A'_{ij} = A'_{rs} = 0$
\item $A'_{is} = A'_{rj} = 1$
\item $A'_{tu} = A_{tu}$ if $t \notin \{i, r\}$ or $u \notin \{j, s\}$.
\end{itemize}
(ignore any points which fall outside of $\Rb$).  In other
words, switch the $i$th and $r$th rows (or equivalently, the $j$th and
$s$th columns) of $\hat{A}$ and restrict to $\Rb$.  See Figure \ref{fundmove}
for an example.

\begin{prop} \label{codone} We have $O_{A'} 
\subset \overline{O_A}$, $\dim(O_A) = \dim(O_{A'}) + 1$,
and all codimension one orbits contained in $\overline{O_A}$ arise this way.
\end{prop} 
\begin{proof}  
The first statement follows from Proposition \ref{Mlorbs}.
To see the second statement, assume first that $(i',j') \in \Rb$.  Then
$\tau(A')$ can be obtained from 
$\tau(A)$ by interchanging rows $i$ and $i'$ and columns $j$ and $j'$,
and then removing the point $(i,j)$.  Thus 
$\#\tau(A') = \#\tau(A) - 1$, as required.  If $(i', j') \notin \Rb$, we
must have either $i' = i+1$ or $j' = j+1$ and a similar argument holds.

For the last part, take $A$, $A'$ so that $O_{A'}$ is a codimension
one subvariety of $O_A$, and let 
$\sigma, \sigma'$ be the corresponding permutations. 
By Lemma \ref{codislen} we have 
$\sigma \le \sigma'$ and 
$\mathop{\rm len}(\sigma) = \mathop{\rm len}(\sigma') - 1$.  A basic result
on reflection groups gives that $\sigma' = \sigma s$, where $s \in S_\infty$
is a reflection.  In other words, $A'$ is obtained from $\hat{A}$ 
by interchanging
two rows and restricting back to $\Rb$.
 
Say the the $i$th and $r$th rows are interchanged, and that $\sigma(i) = j$,
$\sigma(r) = s$.  Then there is no point $(i', j')$ with
$A_{i'j'} = 1$ and $i<i'<r$, $j<j'<s$, since if there were, we would have
 $\mathop{\rm len}(\sigma) < \mathop{\rm len}(\sigma') - 1$.  This gives the
required minimality of $(r,s)$.
\end{proof}

\begin{lemma} \label{nogaps}
For any non-maximal $B_\lambda$ orbit $O \subset M_\lambda$, there is
an orbit $O'$ with $O \subset \overline{O'}$ and $\dim_\C O' =  \dim_\C O + 1$. 
\end{lemma}
\begin{proof} Say $O = O_A$, and 
take a maximal element $(i,j) \in \Rb \setminus \Supp A$.  
Suppose $(i,j'), (i', j) \in \Supp \hat{A}$, and define a matrix
by $A'_{ij} = A'_{i'j'} = 1$, $A'_{ij'} = A'_{i'j} = 0$, $A'_{rs} = A_{rs}$
for all other $(r,s)$.  Then $O_{A'}$ is the required orbit.
\end{proof}

\subsection{Codimension zero orbits} 
Let ${A}_\lambda \in M_\lambda$ be the \zo matrix whose support is
$w(\Pi_{k,l}(\lambda))$, the image of the set of $\lambda$-pairs as defined
in \S1.2, and let $O_\lambda = O_{{A}_\lambda}$.

\begin{prop} \label{codzero} ${O}_\lambda$
is the largest $B_\lambda$-orbit of $M_\lambda$.
\end{prop}
\begin{proof} By Proposition \ref{orbdim}, we need to show that 
$\tau(A) = \Rb$.  Take a point 
$(i,j) \in \Rb$.
Let $(\alpha, \beta) = w^{-1}(i,j)$.
By Lemma \ref{xxx}, we have $\alpha < \beta$.  If $(\alpha, \beta)$ is a 
$\lambda$-pair, we are done.  Otherwise, there are $\lambda$-pairs
$(\alpha, \beta')$ and $(\alpha', \beta)$.  By Lemma \ref{lpnest},
one of these pairs nests inside the other one.  Suppose that
$\alpha < \alpha' < \beta < \beta'$.  Then $(i', j) = w(\alpha', \beta)$
satisfies $A_{i'j} = 1$, and $i' < i$, since $\alpha < \alpha'$.  
Similarly, the other case gives $A_{ij'} = 1$ with $j'<j$.
\end{proof}

The matrices ${A}_\lambda$ were called ``dot configurations" in
\cite{BF}.

Combining Proposition \ref{codzero} with the algorithm of Proposition 
\ref{tatoa}, we obtain:
\begin{lemma} \label{minhole}
If $A \in M_\lambda$ is a \zo matrix, and 
$(i,j)$ is a minimal point in $R_\lambda \setminus \tau(A)$,
then $(i,j)\in w(\Pi_{k,l}(\lambda))$.
\end{lemma}

\subsection{Codimension one orbits} \label{conesec}
Take a $\lambda$-pair $\pi = (\alpha, \beta) \in \Pi_{k,l}(\lambda)$,
and let $\pi' \in  \Pi(\lambda)$ be the parent of $\pi$.
Letting $(i,j) = w(\pi)$, $(r,s) = w(\pi')$, we can use Proposition
\ref{codone} to define a matrix $A_\lambda^\pi$ and a corresponding
codimension one orbit $O_\lambda^\pi$.  Clearly we have
\[\tau(A^\pi_\lambda) = \Rb \setminus \{w(\pi)\}.\]  Thus the
dual orbit $(O_\lambda^\pi)^*$ is represented by the elementary matrix
$E_{w(\pi)}$.

\begin{prop} \label{dualorb}  The correspondence $\pi \mapsto O^\pi_\lambda$
gives a bijection from $\Pi_{k,l}(\lambda)$ to the set of codimension 
one orbits in $M_\lambda$.  

We have $\overline{(O^\pi_\lambda)^*} = \overline{\epsilon^{-1}(X_{\lambda'})}$, 
where $\lambda \ral^\pi \lambda'$.
\end{prop}
\begin{proof} The first statement follows from Proposition \ref{codone}. 
For the second, use Lemma \ref{elemmat}.
\end{proof}

\begin{cor} \label{c1int} $\lambda \lra \lambda'$ if and only if
the corresponding components of $\Lambda$, $\overline{T^*_{X_\lambda}}$ and
$\overline{T^*_{X_{\lambda'}}}$, meet in codimension one.
\end{cor}

\subsection{Codimension two orbits}\label{cod2}

\begin{thm} \label{cod2thm} 
Let $O$ be a codimension two orbit of $M_\lambda$.
There are either one or two codimension one orbits whose closures
contain $O$.  In either case, they have smooth closures at points of $O$.
If there are two, their closures intersect transversely at points of $O$.
\end{thm}
\begin{proof} If $O = O_A$, then $N = A + \C^{\Rb \setminus \tau(A)}$
gives a normal slice to $O$ at $A$.  It is easy to construct
a two-dimensional torus $(\C^*)^2 \subset B_\lambda$ which fixes $A$, 
preserves $N$, and induces the normal crossings stratification on $N$.
Since there must be at least one codimension one orbit $O'$ with 
$O \subset \overline{O'}$ by Lemma \ref{nogaps}, the stratification
on $N$ induced from the orbit stratification of $M_\lambda$ must be
either the normal crossings stratification or the stratification by 
a complete flag.
\end{proof}

We will need a combinatorial parametrization of these codimension
two orbits.  Given $\lambda' > \lambda$ in $\Omega_{k,l}$, let
$O_\lambda^{\lambda'} \subset M_\lambda$ be the dual orbit to the
open orbit in $\epsilon^{-1}(X_{\lambda'})$.

\begin{prop} \label{cod2prop} The map $\lambda' \mapsto O^{\lambda'}_\lambda$ 
defines a one-to-one correspondence between the set of codimension two 
orbits $O \subset M_\lambda$ and the set of $\lambda' \in \Omega_{k,l}$ 
for which there is a diamond (necessarily unique)
$(\lambda \ra \lambda_1 \lra \lambda' \lra \lambda_2 \la \lambda)$
in $\Omega$ where at least one of the $\lambda_m$ is in $\Omega_{k,l}$.  
Given such a diamond and orbit, the orbits $O^{\lambda_m}_\lambda$
are exactly the codimension one orbits whose closures contain $O^{\lambda'}_\lambda$.
\end{prop}
\begin{proof} Take a codimension two orbit $O = O_A$.  Denote the
points in $S = \Rb \setminus \tau(A)$ by $(i_m, j_m)$, $m=1, 2$.  There
are two cases:

\noindent{\bf Case 1:} $i_1 \ne i_2$ and $j_1 \ne j_2$.  Since $A$
must be obtained from ${{A}}_\lambda$ by two applications of 
Proposition \ref{codone}, we have $(i_m, j_m) = w(\pi_m)$ for
$\pi_1, \pi_2 \in \Pi_{k,l}(\lambda)$; neither $\pi_1$ or $\pi_2$ can
be the parent of the other.

The dual orbit $O^*$ contains the \zo matrix with support 
$S$.  We can use Lemma \ref{Mstorbs} to get $\lambda' \in \Omega_{k,l}$
so that $\overline{O^*} = \epsilon^{-1}(\overline{X_{\lambda'}})$ --
if the points in $S$ are comparable in the partial order on $\Rb$,
this is immediate; otherwise it follows from Proposition \ref{tatoa}.
An easy argument along the lines of Lemma \ref{elemmat} shows that
$\lambda_1 \ral^{\pi_2} \lambda'$ and $\lambda_2 \ral^{\pi_1} \lambda'$.

\noindent{\bf Case 2:} The points of $S$ are in the same row or column.
The argument is the same in both cases, so assume WLOG that 
$S = \{(i_1, j_1), (i_2, j_1)\}$ with $i_1 < i_2$.  By Lemma \ref{minhole}
we must have $(i_1, j_1) = w(\pi_1)$ for $\pi_1 \in \Pi_{k,l}(\lambda)$.
Since $(i_2,j_1) \notin \Pi(\lambda)$, we must have 
$w(\pi_2) = (i_2, j_2)$, where $\pi_2 \in \Pi(\lambda)$ is the parent
of $\pi_1$ (note that $\pi_2$ may not be in $\Rb$).

The dual orbit $O^*$ contains $E_{i_2j_1}$, so Lemma \ref{Mstorbs}
gives a partition $\lambda'$ as before.  Lemma \ref{elemmat} now shows that 
$\lambda'$ is the element $\lambda_l$ in the proof of Lemma \ref{diam}.

Finally, the analysis in the proof of Lemma \ref{diam} shows that
all diamonds arise either by case 1 or case 2.
\end{proof}

It is somewhat awkward to index these orbits by diamonds, so we
adopt the following notation, based on the classification of 
diamonds from Lemma\ref{diam}.  If case 1 holds in the proof of the 
previous proposition, we write $O = O^{\pi_1\pi_2}_\lambda$.
In case 2, put $O = O^{\pi_1,l}_\lambda$ or $O^{\pi_1,r}_\lambda$
if the points of $S$ are in the same column or row, respectively.

\begin{figure}
\begin{center}
\leavevmode
\hbox{
\epsfxsize=3.5in
\epsffile{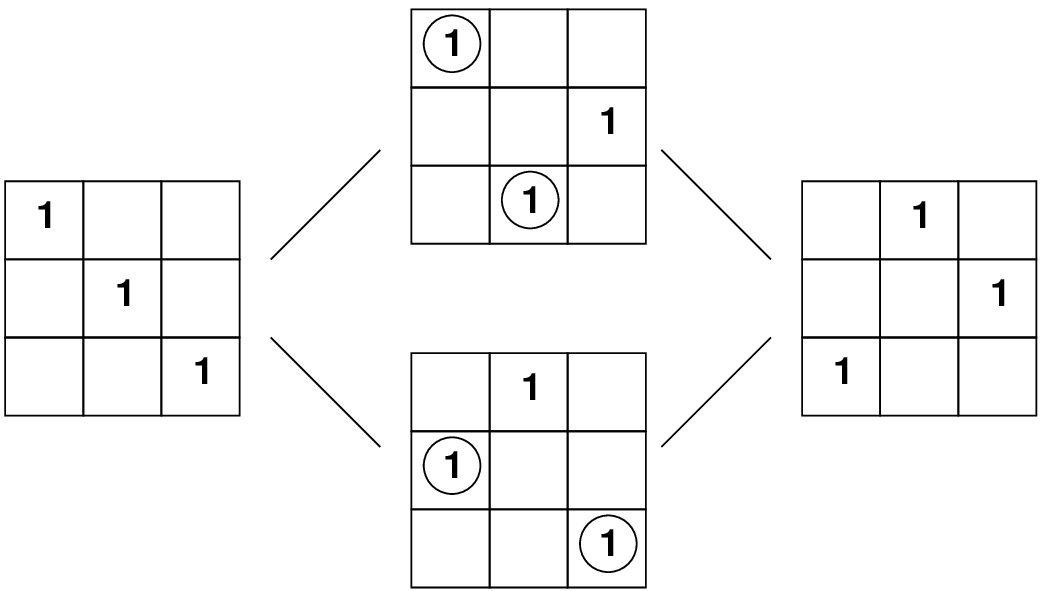}}
\end{center}
\caption{}
\label{typeII}
\end{figure}

Figure \ref{typeII} illustrates case 2 with an example: let $k = l = 3$, 
$\lambda = \emptyset$, and let $\pi_1 = (-\half, \half)$, so $\pi_2 = 
(-\frac{3}{2}, \frac{3}{2})$.  The rightmost
matrix represents the stratum $O^{\pi_1,l}_\lambda$.
The corresponding diamond is $\{\emptyset, (1), (2, 1), (1, 1)\}$.

\subsection{Fundamental groups of orbits} \label{pi1} 
\newcommand{\Tb}{{\overline{T}}}

The stabilizer $(B_\lambda)_A$ of a \zo matrix $A \in M_\lambda$ is
given by intersecting $B_\lambda$ with a linear subspace, and is thus
connected.
As a result, the fundamental group $\pi_1(O_A,A)$
is a quotient of $\pi_1(B_\lambda, 1)$; 
we can generate $\pi_1(O_A)$ by acting on $A$
by loops in $\pi_1(B_\lambda) \cong \pi_1(T)$, 
where $T \subset B_\lambda$ is the torus of diagonal matrices.  
Since these groups are abelian, we can safely ignore basepoints.

For the sake of convenient notation,
it is easier to work with the larger torus $\overline{T} \cong (\C^*)^{k+l}$
of all diagonal matrices.  Define a collection of generators for
$\pi_1(\overline{T})$, indexed by elements 
$\alpha \in \Hi \cap (-k, l)$: let
\[\gamma_\alpha(e^{i\theta}) = (1, \dots, e^{i\theta}, \dots, 1),\]
where the $e^{i\theta}$ is in the $(\alpha + k + \half)$th place.  

The element $\gamma_\alpha(z)\in \overline{T}$ acts on $M_\lambda$ 
as follows:  if $\phi_\lambda(\alpha) = -1$, multiply the row numbered
$w_1(\alpha)$ by $z$.  If $\phi_\lambda(\alpha) = +1$, multiply column
number $w_2(\alpha)$ by $z^{-1}$ ($w_1$ and $w_2$ are
the component functions of $w$; see \S\ref{ncc}).  
The action on $M_\lambda^*$ is the same, with $z$ and $z^{-1}$ interchanged. 

It follows that if $O = O_A$ is a $B_\lambda$-orbit in $M_\lambda$,
then $\pi_1(O) \cong \Z^{\#\Supp A}$.  In particular, if we let
$O = {O}_\lambda$, the rank of 
$\pi_1(O)$ is $\#\Pi_{k,l}(\lambda)$. We will abuse notation
and use the same symbol $\gamma_\alpha$ to denote a loop in $\pi_1(\Tb)$ 
and its image in an orbit $O$.  

\begin{prop} \label{lsq} The fundamental group $\pi_1(O_\lambda)$
is the abelian group generated by the $\gamma_\alpha$,  $\alpha \in \Hi$,
with relations $\gamma_\alpha\gamma_\beta = 1$ if $(\alpha, \beta)\in \Pi(\lambda)$, and 
$\gamma_\alpha = 1$ if $\alpha \notin [-k,l]$.
\end{prop}
Thus a local system on $O_\lambda$ can be described as a representation
of the group algebra $\C\pi_1(O_\lambda)$, 
which has commuting generators
$t_\alpha$, $\alpha \in \Hi$ satisfying the relation (2) 
from the quiver algebra $A_{k,l}$ restricted
to $V_\lambda$. 

Finally, given $O = {O}_\lambda$, and a codimension $1$ 
orbit $O' = O_\lambda^\pi$ in $M_\lambda$, we wish
to describe the class in $\pi_1(O)$ of a loop around a point of $O'$,
following the holomorphic orientation of the normal slice to $O'$.
Let $\pi = (\alpha, \beta)$, and let $\pi' = (\alpha', \beta')$ be the
parent of $\pi$.
A normal slice to $O'$ through the \zo matrix $A' \in O'$
is given by $N = A' + \C^{w(\pi)}$.  Since $\gamma_\alpha(z)
\gamma_{\beta'}(z)$ multiplies the entry at $w(\pi)$ by $z$ and
fixes $A'$, this gives the required loop.  This calculation will be
important in \S4 when we explain relation (3) from the quiver
algebra $A_{k,l}$.

\subsection{Modifications for type D}  The analysis of the geometry
of $X = X^s_k$ proceeds very similarly to the preceding discussion.
We will only indicate the places where the arguments must be changed.
Keeping the same definition of $W = W_\lambda$, we have $W_\lambda \in X^s_k$
if and only if $\lambda \in \Omega^s_k$.  It is the unique point in
$X^s_\lambda$ which is fixed by the torus $T\subset SO(n, \C)$ of diagonal
matrices.

\subsubsection{Coordinates} 
The symmetric bilinear form associated to the quadratic form $Q$ 
gives an identification 
of $\C^n$ with $(\C^n)^*$, which in turn 
gives an identification $\C^n/W \cong W^*$
when $W$ is isotropic. 
Using the identification $T_WX_{k,k} = \Hom(W, \C^n/W) \cong \Hom(W, W^*)$,
we have 
\[T^{}_WX^s_k = \{\phi \in \Hom(W, W^*) \mid \phi + \phi^t = 0\},\]
so using the standard basis of $W$ and the dual basis of $W^*$,
tangent vectors are represented by skew-symmetric matrices.  

Let $R =  \{1, \dots, k\}\times\{1, \dots, k\} \setminus 
\{(i,i) \mid 1\le i\le k\}$.
For a subset $S\subset R$ we let 
$S^t = \{(j,i) \mid (i,j) \in S\}$.  If $S = S^t$, 
define $\C_s^S$ to be the set of skew-symmetric matrices supported on $S$.
We then have $T_WX^s_\lambda = \C_s^{\Delta(\lambda)}$, and so if $M_\lambda$
is the fiber $M_\lambda$ of the conormal variety
$\Lambda$ over $W_\lambda$ we get an identification
$M_\lambda \cong \C_s^\Rb$, where $\Rb = R \setminus \Delta(\lambda)$.

Using the standard pairing, the space of skew-symmetric matrices is dual 
to itself, so we also have an identification $M^*_\lambda \cong \C_s^\Rb$.

\subsubsection{Actions and orbits}
If $B_\lambda = B_{W_\lambda}$, then the action of $B_\lambda$ on $M_\lambda$ is 
generated by ``symmetric elementary
operations", in which the same operation is performed to both the rows and 
the columns of a matrix, and rows or columns are only added to higher-numbered
rows or columns.  The dual action on $M^*_\lambda$ is given by truncated versions
of these operations, where rows or columns are only added to lower-numbered
ones, and anything appearing in $\Delta(\lambda)$ is discarded.

Call a (skew-symmetric) matrix in $M_\lambda$ or
$M^*_\lambda$ a \zo matrix if all of its entries above the diagonal 
are $0$ or $1$ and each row or column has at most one nonzero
element.  Then everything between Proposition \ref{Mlorbs} and 
 Proposition \ref{hjk} holds word-for-word 
in type D, including the definition of the embedding 
$\epsilon\colon M^*_\lambda \to X^s$.  In particular $A \mapsto O_A = B_\lambda A$
defines a one-to-one correspondence
between $B_\lambda$ orbits in $M_\lambda$ (or $M^*_\lambda$) and \zo matrices, as
 before.

Lemma \ref{Mstorbs} must be modified in the type D case:
\begin{lemma} \label{DMstorbs} Take a \zo matrix $A$, considered as an
element of $M^*_\lambda$.  Then $B_\lambda A$ is the 
(unique) open orbit in some $\epsilon^{-1}(X_{\lambda'})$ if and only if
for each $(i, j)$, $(r, s) \in \Supp A$ with $i<r$, $j>s$
{\bf and $\mathbf{i {\pmb \ne} s}$} we have $(i, s) \in \Delta(\lambda)$.
\end{lemma}

Given a point $(i,j) \in \Rb$, we define an ``elementary matrix"
$E^s_{ij} = E_{ij} - E_{ji} \in M^*_\lambda$.  Lemma \ref{DMstorbs}
gives $\lambda' \in \Omega^s_k$ so that $B_\lambda E^s_{ij}$ is the 
unique open orbit in $\epsilon^{-1}(X_{\lambda'})$.  Let $(\alpha, \beta) =
w^{-1}(i,j)$.  Just as in Lemma \ref{elemmat}, we have:

\begin{lemma} $\phi_\lambda(\gamma) = \phi_{\lambda'}(\gamma)$ if and only
if $\gamma \notin \{\pm\alpha, \pm\beta\}$.
\end{lemma}

Define $\tau(A)$ for a \zo matrix $A$ exactly as in \S\ref{orbstrs} 
(remember that the diagonal has been removed from $R$).

\begin{prop} If $A$ is a \zo matrix, we have 
\[T_A(O_A) = \C_s^{\tau(A)}.\]
Thus $\dim O_A = \half \#\tau(A)$.  The algorithm of Proposition
\ref{tatoa} works verbatim, except that $-1$'s are placed 
below the diagonal instead of $1$s.
\end{prop}

As before, we can extend a \zo matrix $A$ to a matrix $\hat{A}$ 
on all of $\N \times \N$,
by applying the algorithm of Proposition \ref{tatoa} to the set
\[\tau(A) \cup \{(i, j) \mid i\ne j\;\text{and}\;(i>k \;\text{or}\; j >k)\}.\]

We can again describe all orbits $O\subset \overline{O_A}$ with 
$\dim O = \dim O_A - 1$, but the procedure is slightly more complicated than 
in Proposition \ref{codone}.  Begin as before with a point $(i,j) \in \Supp A$,
and assume that $i < j$.  Let $(i',j')$ be a minimal 
point in $\Supp \hat{A}$ so that $(i,j) < (i',j')$. If $j < i' < j'$, then
switch $i'$ and $j'$.  In geometric terms, this ensures 
that no corner of the rectangle
$\rho$ with corners at $(i,j)$ and $(i',j')$ lies in the reflected rectangle
${\rho}^t$.
Given these two points, define a \zo matrix $A'$ by switching the $i$th and
$i'$th columns and the $i$th and $i'$th rows, multiplying any squares by
$-1$ that are needed to make it a \zo matrix. 

\begin{prop} \label{codones} We have $O_{A'} 
\subset \overline{O_A}$, $\dim(O_A) = \dim(O_{A'}) + 1$,
and all codimension one orbits contained in $\overline{O_A}$ arise this way.
\end{prop}

\begin{figure}
\begin{center}
\leavevmode
\hbox{
\epsfxsize=3.7in
\epsffile{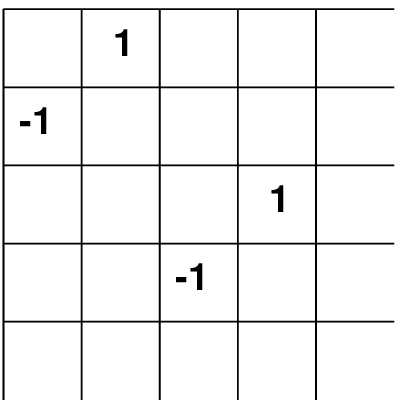}}
\end{center}
\caption{}
\label{xyzzy}
\end{figure}

As an example of the extra step in this procedure, consider
the left matrix in Figure \ref{xyzzy}, representing the largest orbit in 
$M_\emptyset$.  Let $(i,j) = (1, 2)$ and $(i',j') = (3, 4)$.  Then switching
the first and third rows and columns (the center matrix) gives a codimension
two orbit, while switching the first and fourth columns and rows gives the 
right matrix, which represents a codimension one orbit.

\subsubsection{Codimension zero and one orbits}

Let ${{A}}_\lambda$ be the \zo matrix that lies in the open orbit of
$M_\lambda$, and let $O_\lambda = O_{{{A}}_\lambda}$.

\begin{prop} \label{Dtopmat} We have $\Supp {{A}}_\lambda = w(\Pi^s_k(\lambda)) \cup
[w(\Pi^s_k(\lambda))]^t$.
\end{prop}
\begin{proof} Recalling the definition of $\Pi^s_k(\lambda)$, we see that
points in $w(\Pi^s_k(\lambda))$ are of two types: first, all points
$(i,j) \in w(\Pi_{k,k}(\lambda))$ with $i<j$, and second, points $(i,j)$
where $(i,i)$, $(j,j)\in \Pi_{k,k}(\lambda)$ run over all
 points in $w(\Pi_{k,k}(\lambda))$ which lie on the diagonal,
taking these points in pairs 
(the first and second, then the third and fourth, and so on).

It is clear from this description that there is a \zo matrix $A$ with
the required support, and that $\tau(A) = \Rb$.
\end{proof} 

There is a one-to-one correspondence
$\pi \mapsto O^\pi_\lambda$ between $\Pi^s_k(\lambda)$ and the set of
codimension one orbits in $M_\lambda$, defined as follows.  
Let $\pi'$ be the parent of $\pi$ in $\Pi^s(\lambda)$, and set 
$(i,j) = w(\pi)$, $(i',j') = w(\pi')$.  Apply the transformation of
Proposition \ref{codones} to the matrix ${A}_\lambda$ using these 
two points; call the resulting matrix $A^\pi_\lambda$, and set
$O^\pi_\lambda = O_{A^\pi_\lambda}$.

The description of the dual orbit to $O^\pi_\lambda$ given in
Proposition \ref{dualorb} works verbatim.

\subsubsection{Codimension two orbits}
Theorem \ref{cod2thm} remains true in the type D case, with essentially 
the same proof.  Proposition \ref{cod2prop} must be changed slightly, however.
Every diamond $(\lambda \ra \lambda_1 \lra \lambda' \lra \lambda_2 \la \lambda)$
gives rise to a codimension two orbit $O = O^{\lambda'}_\lambda$, but there
is another kind of codimension two orbit.  If $\pi \in \Pi^s_k(\lambda)$,
and the parent of $\pi$ is not in $\Pi(\lambda)$,  the two points in 
$\Supp A^{\pi}_\lambda \setminus \Supp {{A}}_\lambda$ that
lie above the diagonal are parent and child, rather than
siblings as in Proposition \ref{cod2prop}.  This is because the modified
rule in Proposition \ref{codones} must be applied.

Applying Proposition \ref{codones} to these 
two points gives, as before, two codimension two orbits (call them 
$O^{\pi,p}_\lambda$ and $O^{\pi,c}_\lambda$).  
The first orbit behaves as in the type A case, but the orbit 
$O^{\pi,c}$ does not correspond to a diamond -- $O^\pi_\lambda$ is
the only codimension one orbit whose closure contains $O^{\pi,c}_\lambda$,
even if $k$ is large. The orbits $O^{\pi,c}_\lambda$ are responsible 
for the extra clause in the quiver relation 4s.

For instance, take $\lambda = \emptyset$. If $\pi_1 = (-\half, \frac{3}{2})$,
$\pi_2 = (-\frac{5}{2}, \frac{7}{2})$, then $O^{\pi_1\pi_2,c}_\lambda$
corresponds to the sequence of partitions $\emptyset \to (2, 2) \to (3, 2, 1)$;
it is easy to check that this cannot be completed to a diamond.

\subsubsection{Fundamental groups} \label{Dpi1} 
As in the type A case the
fundamental group of an orbit $O_A$ is a quotient of $\pi_1(T)$, where
$T\subset B_\lambda$ is the torus of diagonal matrices.  
In terms of the loops $\gamma_\alpha$ described in \S\ref{pi1}, generators
for $\pi_1(T)$ are given by
 $\bar{\gamma}_\alpha = \gamma_\alpha\gamma^{-1}_{-\alpha}$, for
$-k< \alpha < k$.
  
\begin{prop} \label{Dlsq} $\pi_1(O_\lambda)$
is the abelian group generated by the $\bar\gamma_\alpha$,  $\alpha \in \Hi$,
with relations $\bar\gamma_\alpha\bar\gamma_{-\alpha} = 1$, 
$\gamma_\alpha\gamma_\beta = 1$ if $(\alpha, \beta)\in \Pi^s(\lambda)$, and 
$\gamma_\alpha = 1$ if $\alpha \notin [-k,k]$.
\end{prop}
In other words, the $\bar\gamma_\alpha$ satisfy the relations 
(2s) for the quiver algebra $\cA^s_k$ (\S1.7).

The relation (3s) in the quiver describes
the class of a loop in the open orbit ${O}_\lambda$
around a point of the codimension one orbit $O^\pi_\lambda$ 
in terms of the generators $\bar{\gamma}_\alpha$.  Recalling
the construction of $A_\lambda^\pi$, we see that 
$\bar{\gamma}_\alpha\bar{\gamma}_\zeta$
gives the required loop, where $\pi = (\alpha, \beta)$, and $\zeta$ is
as described in (3s).

\section{Microlocal perverse sheaves}
In this section we collect some facts about microlocal perverse sheaves that we will use 
in our calculation, along with their local description along conormal varieties that have
normal crossings singularities through codimension two.

\subsection{Stacks}
Our calculation of perverse sheaves rests on the observation that 
perverse sheaves on a complex manifold $X$ can be seen as microlocal
objects, i.e., objects which are locally defined on $T^*X$.  To say this
more precisely, we need the language of stacks.  We give a sketch of 
a simplified definition which is sufficient for our application.

A stack on a topological space 
$Y$ can be thought of intuitively as a ``sheaf of categories'' on $Y$.
Formally, it consists of a category $\cS(U)$ for every open $U \subset Y$,
together with a restriction functor $|^U_V\colon \cS(U) \to \cS(V)$ for 
any open subset $V\subset U$  (when no confusion will arise, we write simply $|_V$).  
They are required to be compatible in the sense that if $W \subset V \subset U$, then
\[|^V_W \circ |^U_V = |^U_W.\]

\begin{rem}
Strictly speaking one cannot say that two functors are equal, so
a formally correct definition would include 
natural isomorphisms relating these two functors, which then must satisfy a sort
of cocycle identity.  However, for our stacks $\cS(U)$ can be considered
as a subcategory of the category of modules over $\cR(U)$ for some sheaf of rings
$\cR$ on $Y$, and  
in this case it makes sense to say that the functors are equal.
A correct general definition can be found in the Appendix to \cite{GMV}.
\end{rem}

In order for $\cS$ to be a stack, objects and morphisms must be locally 
defined.  Let $U\subset Y$ be an open set, and let $\cU = \{U_\alpha\}_{\alpha \in A}$ be a cover of 
$U$ by open sets.  
We put $U_{\alpha\beta} = U_\alpha \cap U_\beta$, and $U_{\alpha\beta\gamma} = U_\alpha \cap U_\beta
\cap U_\gamma$.  Define
a category $\cS(\cU)$ by letting objects be collections $\{F_\alpha\}_{\alpha\in A}$ 
where $F_\alpha$ is an object of $\cS(U_\alpha)$, together with isomorphisms
\[\phi_{\alpha\beta}\colon F_\beta|_{U_{\alpha\beta}} \to F_\alpha|_{U_{\alpha\beta}},\,
\alpha,\beta\in A.\]
They are required to satisfy
\[\phi_{\alpha\beta}|_{U_{\alpha\beta\gamma}} \circ 
\phi_{\beta\gamma}|_{U_{\alpha\beta\gamma}} =
\phi_{\alpha\gamma}|_{U_{\alpha\beta\gamma}},\,
\alpha,\beta,\gamma\in A.\]
A morphism $f\colon (\{F_\alpha\}, \phi_{\alpha\beta}) \to (\{G_\alpha\}, \psi_{\alpha\beta})$
is a collection of morphisms $f_\alpha\colon F_\alpha \to G_\alpha$ satisfying \[\psi_{\alpha\beta} \circ 
f_\beta|_{U_{\alpha\beta}} 
= f_\alpha|_{U_{\alpha\beta}} \circ \phi_{\alpha\beta}\]
for all $\alpha, \beta \in A$.

\begin{defn} We say that $\cS$ is a stack if for all $U$ and $\cU$
as above, the obvious functor $\cS(U) \to \cS(\cU)$ given by letting $F_\alpha = F|^U_{U_\alpha}$
is an isomorphism of categories.
\end{defn}

A {\em substack} 
$\cS'$ of $\cS$ is a stack for which each $\cS'(U)$ is a full subcategory
of $\cS(U)$.  We call $\cS$ a stack of abelian categories if each $\cS(U)$ is 
abelian, and the restriction functors are exact.  In that case, given an object
$A\in \cS(U)$, we define its support $\Supp A\subset U$ to be 
the smallest closed set for which 
$A$ restricted to $U \setminus \Supp A$ is zero.

\begin{rem} The categories $\cS h(X)$ of sheaves and $\cP(X)$ of perverse sheaves
on $X$ both form stacks of abelian categories on $X$.  
The derived category $D^b(X)$, however, is {\em not} 
generally a stack.  For instance, if $X = \C\PP^1$, there is a nonzero morphism
$\Q_X \to \Q_X[2]$ whose restriction to any smaller open set is zero.
\end{rem}

\subsection{The stack $\cE$}

Microlocal perverse sheaves form a stack $\cE$ of abelian categories
on the cotangent bundle $T^*X$.  It has
been constructed topologically in \cite{GMV}, or one can use the 
stack of holonomic regular singularities
 $\cE_X$-modules, via the Riemann-Hilbert
correspondence \cite{A, K2,K,KK}.  
It has the following properties:
\begin{itemize}
\item
If $U\subset X$ is open, there is an equivalence of categories
$\cE(T^*U) \cong \cP(U)$.  
\item The support of 
any object in $\cE(U)$ is a locally $\C^*$-invariant Lagrangian variety in $U$.
\item 
The support
of an object $A \in \cE(T^*X)$ is equal to its micro-support
$SS(A)$, considering it
as an object in $\cP(X)$.  
\end{itemize}

Given a conical Lagrangian variety $\Lambda\subset T^*X$,
we let $\cE_\Lambda$ be the substack of objects supported on $\Lambda$.
Then we have an equivalence of categories
\[\cE_\Lambda(T^*X) \cong \cP_\Lambda(X).\]
Thus, if $\Lambda$ is the conormal variety to a stratification $\cS$ of $X$,
then $\cE_\Lambda(T^*X)$ is equivalent to the category of
 $\cS$-constructible perverse sheaves on $X$.

The following theorem is proved in \cite{GMV} (also see \cite{K}):
\begin{thm} \label{GMVcod1}If $U\subset T^*X$ is a Zariski open set and
$\dim_\C (\Lambda\cap (T^*X \setminus U)) \le 
\dim_\C \Lambda-2$, then
the restriction functor \[\cP_\Lambda(X) \to \cE_\Lambda(U)\] is a full embedding 
of categories.
\end{thm}

In other words, microlocal perverse sheaves through codimension one
have enough information to completely describe perverse sheaves and
maps between them, but
there may be further relations which must be satisfied in order to extend
an object from $U$ to all of $T^*X$.

Kashiwara has conjectured that microlocal perverse sheaves through 
codimension two are equivalent to global perverse sheaves.  We
need the following weaker result:

\begin{thm} \label{weakKash} If $U\subset T^*X$ is a Zariski open set and
 $\dim_\C (\Lambda\cap (T^*X \setminus U)) \le 
\dim_\C \Lambda-3$, then 
\[\Ext^1_{\cP_\Lambda(X)}(A, B) = \Ext^1_{\cE_\Lambda(U)}(A|_U,B|_U)\]
for any two objects $A, B$ in $\cP_\Lambda(X)$.
\end{thm}
\begin{proof} This follows directly from \cite{KK}, Theorem 1.2.2.  
Also see \cite{K}, Theorem 10.4.1.
\end{proof}

\subsection{Microlocal Fourier transform}

Given a complex vector bundle $E$, the Fourier transform is a functor
$F\colon D^b_c(E) \to D^b_c(E^*)$, where $D^b_c(E)$ denotes objects
of th derived category which are 
constructible with respect to some analytic $\C^*$-conic stratification.  We will
shift our Fourier transform by $\dim_\C E$
from the definition in \cite{KS}, 
so that it sends perverse sheaves to perverse sheaves.
Let $\Theta\colon T^*E \cong T^*E^*$ be the natural identification
(\cite{K}, Proposition 5.5.1).

\begin{prop} \label{mlft} 
$F$ preserves the micro-support: $SS(F{\mb P}) = \Theta(SS({\mb P}))$.
Given $U\subset T^*E$ open and the conormal variety $\Lambda$
to a conical stratification of $E$, $F$ 
induces an isomorphism of categories 
\[\cE_\Lambda(U) \to \cE_{\Theta(\Lambda)}(\Theta(U)).\]
\end{prop}

\subsection{Monodromic perverse sheaves}
We need to describe carefully the categories of perverse
sheaves which give the local descriptions we will glue together.
We begin with a model for the codimension one data.

Let $L$ be a complex vector bundle over a smooth connected base $B$. 
Let $\cP_{mon}(L)$ be the 
category of {\em monodromic\/} perverse sheaves, i.e.\ sheaves which are constructible
with respect to the stratification $(Z, \wt{L})$, where $Z$ is the zero 
section, and $\wt{L} = L \setminus Z$.

We briefly describe the standard model for $\cP_{mon}(L)$.
Let $L'$ be the dual vector bundle and let $\wt{L}'$ be the bundle minus the
zero section.  We can define a isomorphism of $\C^*$ bundles 
$\omega\colon \wt{L} \to \wt{L}'$
by requiring that $\langle \omega(x), x\rangle = 1$ for all $x$.  Given
local systems $\cL$, $\cL'$ on $\wt{L}$ and $\wt{L}'$ respectively,
let $\mu\colon \cL\to \cL$ and
$\mu'\colon \cL' \to \cL'$ be the automorphisms given by following the
monodromies around loops given by multiplication by
$e^{i\theta}$ and $e^{-i\theta}$, respectively.  

\begin{lemma}[\cite{MV}, \cite{V}]\label{Pmon}
There is a natural equivalence of categories  $R_{mon}\colon \cP_{mon}(L) \to \cC$, where 
$\cC$ is the category
whose objects are local systems $\cL$ and $\cL'$ on $\wt{L}$ and $\wt{L}'$ respectively,
together with morphisms of local systems
\[p\colon \omega_*\cL \to \cL',\;\; q\colon \cL' \to \omega_*\cL\]
satisfying 
\[pq+1 = \mu',\;\text{and}\; qp+1 = \omega_*\mu.\]
\end{lemma} If ${\mb P} \in \cP_{mon}(L)$, then
these local systems are described concretely by
\[\cL = \cL({\mb P}) = {\mb P}|^{}_{\wt{L}};\]
\[ \cL' = \cL'({\mb P}) = (F{\mb P})|^{}_{\wt{L}'},\]  
where $F\colon \cP_{mon}(L) \to \cP_{mon}(L')$ is the Fourier transform.
The maps $p$, $q$ are described in \cite{MV}, \cite{V}.
Note that they depend on the choice of a fixed orientation of $L$.  

If $L$ is oriented, let $L'$ have the orientation induced using
$\omega$, so if $L$ has the holomorphic orientation, $L'$ has the
antiholomorphic one.  Then Lemma \ref{Pmon} gives an equivalence
$R'_{mon}\colon \cP_{mon}(L')\to \cC'$, where $\cC'$ is the category
defined above taking $L'$ instead of $L$.  If we define a Fourier
transform $F_\cC\colon \cC \to \cC'$ by
\[F_\cC(\cL, \cL', p, q) = (\cL', \cL, \omega^*q, \omega^*p),\]
then there is a natural equivalence $F_\cC \circ R_{mon} \cong R'_{mon} \circ F$.

\subsection{Normal crossings perverse sheaves}\label{ncpss}

Let $L_1$ and $L_2$ be complex line bundles over a smooth, connected base $B$, and let
$E = L_1 \oplus L_2$.  Give $E$ the ``normal crossings'' stratification $\{S_{00}, S_{01}, S_{10}, S_{11}\}$,
where $S_{00} = L_1 \cap L_2$ is the zero section,
$S_{11} = E\setminus(L_1 \cup L_2)$, $S_{10} = L_1 \setminus S_{00}$, and
$S_{01} = L_2 \setminus S_{00}$.  Let $\Lambda = \Lambda_{nc}$ 
be the corresponding conormal variety, and let ${\wt \Lambda}_{ij}$ be the component of 
the smooth locus of $\Lambda$ lying over $S_{ij}$.

Set $U^1 = T^*E \setminus \{(x, 0) \in T^*E \mid x \in S_{00}\}$, so
$U^1$ contains the codimension zero and one parts of $\Lambda$.  
By Theorem \ref{GMVcod1}, $\cP_\Lambda(E)$ embeds as a full subcategory
of $\cE_\Lambda(U^1)$.  We wish to describe the category 
$\cE_\Lambda(U^1)$ and the additional 
relations that objects must satisfy in order to extend to $\cP_\Lambda(E) = \cE_\Lambda(T^*E)$.
 
First we reduce to the case where $B$ is a point. Take $b\in B$, and let
$i_b\colon E_b \to E$ denote the inclusion of the fiber over $b$.  Give
$E_b$ the induced stratification, and define $\Lambda_b, U^1_b \subset T^*E_b$
as above.  Then the restriction functor 
$i^*[-\dim_\C B]\colon D^b(E) \to D^b(E_b)$ takes perverse sheaves to 
perverse sheaves, and microlocalizes to give a functor
\[i_\mu^*\colon \cE_\Lambda(U^1) \to \cE_{\Lambda_b}(U^1_b).\]

\begin{lemma} \label{onec2pt} An object ${\mb P} \in \cE_\Lambda(U^1)$ 
extends to $\cP_\Lambda(E)$ if
and only if $i_\mu^*{\mb P}$ extends to $\cP_{\Lambda_b}(U_b)$.
\end{lemma}
\begin{proof}
This is easy to see when $B$ is contractible, since the bundles $L_1$ and 
$L_2$ can be trivialized.  For the general case, cover $B$ by contractible open sets 
containing $b$.  Since extensions from $\cE_\Lambda(U^1)$ to $\cE_\Lambda(T^*E)$ are unique
and canonical if they exist, by Theorem \ref{GMVcod1}, the extensions over these open 
sets glue together.
\end{proof}

Now restrict to the case where $E = \{b\}$ is a point, so $L_1$ and $L_2$ are
just 1-dimensional vector spaces.  To avoid naming every map in the 
quiver diagram, we use the convention that 
$(v_1\mid  \dots \mid v_r)$ denotes the 
composition of the maps along the path with nodes labelled $v_1, \dots, v_r$.

\begin{prop} \label{ncps} There is an equivalence of categories $R_{nc}\colon\cE_\Lambda(U^1) \to \cC_{nc}$,
where $\cC_{nc}$ is the category of representations of the quiver
\begin{equation} \label{ncquiver}
\xymatrix{
V_{00} \ar@<.5ex>[r] \ar@<.5ex>[d] & V_{01} \ar@<.5ex>[l]\ar@<.5ex>[d]\\
V_{10} \ar@<.5ex>[u] \ar@<.5ex>[r] & V_{11} \ar@<.5ex>[l]\ar@<.5ex>[u]}
\end{equation}
which satisfy 
\[(00 \mid 01 \mid 00 \mid 10) = (00 \mid 10 \mid 11 \mid 10),\,\text{and}\;
(00 \mid 01 \mid 00) \;\text{invertible}\] as well as all relations obtained from these
by applying the symmetry group of the square.

The category $\cP_\Lambda(E)$ sits inside $\cE_\Lambda(U^1)$ as quivers satisfying
the additional relations of the form \[(00 \mid 10 \mid 11) = (00 \mid 01 \mid 11),\]
again applying the symmetries of the square.
\end{prop}

\begin{proof} Fix orientations of $L_1$ and $L_2$, and give the
dual bundles $L_1'$
and $L_2'$ the compatible orientations, as in the last section. 
Using the description of $\cP_{mon}$ from the last
section twice, we see 
that an object ${\mb P} \in \cP_\Lambda(E \setminus S_{00})$ 
is described by local systems
$\cL_{ij}$ on ${\wt \Lambda}_{ij}$ for $ij = 11, 01, 10$,
together with maps (making identifications between these spaces as in Lemma
\ref{Pmon})
$\cL_{01} \leftrightarrows \cL_{11} \leftrightarrows \cL_{10}$;
these maps must satisfy relations as in Lemma \ref{Pmon}.
In particular, the monodromies of $\cL_{11}$ are completely determined 
by what these maps do on stalks of the $\cL_{ij}$.

Applying the Fourier transform and 
using Proposition \ref{mlft}, we see that
objects in $\cE_\Lambda(T^*E \setminus \overline{{\wt \Lambda}_{11}})$
are given by local systems and maps $\tilde\cL_{01} 
\leftrightarrows \tilde\cL_{00} \leftrightarrows \tilde\cL_{10}$, 
with relations as before.  

We can glue these two descriptions using the fact that $\cE$ is a stack
to get a description of the category
$\cE_\Lambda(U^1)$; the key fact is that there are natural
identifications $\tilde\cL_{ij}({\mb P}) \cong \cL_{ij}({\mb P})$, 
for $ij = 01$ and $10$.
This follows easily from existence of natural isomorphisms (\cite{BG},
Proposition 2.4)
$F \cong F_1F_2 \cong F_2F_1$, where $F_1$ denotes the Fourier transform
functors $\cP_\Lambda(L_1 \oplus K) 
\leftrightarrows \cP_\Lambda(L_1' \oplus K)$, $K = L_2$ or $L_2'$,
 and similarly $F_2$ is the Fourier transform in the second coordinate.

Now fix a basepoint $p_{11} \in {\wt \Lambda}_{11}$; let 
$p_{ij} \in {\wt \Lambda}_{ij}$ be the corresponding points in the other
smooth components.
The vector space $V_{ij}$ in the quiver (\ref{ncquiver}) 
is the stalk of $\cL_{ij}$ at $p_{ij}$, and the map $V_{ij} \to V_{kl}$
is induced from the map $\cL_{ij} \to \cL_{kl}$.
The maps on stalks determine the monodromies of all four
local systems, and the quiver relations arise because 
the $\cL_{ij}$ are local systems, and the maps are maps of local systems.

The description of the category $\cP_\Lambda(E)$ is standard; see
\cite{V}, \cite{GGM}.
\end{proof}

We will need a more microlocal description of the orientation choices 
we made.  Instead of choosing orientations of the $L_i$, we can give each 
component of $\Lambda$ a normal orientation along each codimension one
intersection with another component.  This gives eight orientations; 
there are four compatibility relations required to apply 
Lemma \ref{Pmon}, and there are two further relations needed to 
make Proposition \ref{ncps} work.  For instance, normal
orientations of $\Lambda_{00}$ along $\Lambda_{00} \cap \Lambda_{01}$
and of $\Lambda_{10}$ along
$\Lambda_{10} \cap \Lambda_{11}$ both give orientations of $L_2'$; these
must agree.

As a consequence 
of the proof of Proposition \ref{ncps}, we can also describe the action 
of the Fourier transform on the quiver (\ref{ncquiver}).  If we identify $L_i \cong L_i'$ 
antiholomorphically, we
get functors $F_1$, $F_2$, $F\colon \cE_\Lambda(U^1) \to \cE_\Lambda(U^1)$;
they are naturally equivalent under $R_{nc}$ to the functors given by reflecting the diagram
(\ref{ncquiver}) in a horizontal line, reflecting in a 
vertical line, and rotating by 180$^\circ$, respectively.

\section{Proof of the main theorems}
We now apply the results on microlocal perverse sheaves of \S3 to the geometry of
\S2 to prove our main theorems.

\subsection{Reducing the theorems to microlocal geometry}
Let $X$ be one of the spaces $X_{k,l}$ or $X^s_k$, taken with the
Schubert stratification and corresponding conormal variety $\Lambda$. 
By the results of \S\ref{mlgeom}, $\Lambda$ has finitely many orbits under the action of $B$.  
Recall that there is a one-to-one correspondence between $B$-orbits of
$\Lambda$ and the union of all $B_\lambda$-orbits of $M_\lambda$ over all
$\lambda$.  Using the notation for the codimension $0$, $1$, and $2$ orbits of
$M_\lambda$ introduced in \S2, we 
let $\cO^\pi_\lambda \subset \Lambda$ be the orbit 
corresponding under Proposition  2.2.1 to 
$O^\pi_\lambda \subset M_\lambda$, and so on.  In particular,
$\cO_\lambda = \wt\Lambda_\lambda$ is the smooth component of $\Lambda$ lying
over $X_\lambda$.

Let $U^d$ be the union of $T^*X \setminus \Lambda$ and all 
$B$-orbits in $\Lambda$ of codimension at most $d$.  The following
result explains how the quiver categories from \S1.3 and \S1.7 arise.

First we show how to deduce our theorems from the following statement.

\begin{thm} \label{abcd} In the case $X = X_{k,l}$ there is an equivalence of categories
\[R\colon \cE_\Lambda(U^2) \to \cA_{k,l}.\]  If $X = X^s_k$, there is an equivalence of categories
\[R\colon \cE_\Lambda(U^2) \to \cA^s_k.\]  In both cases, if $\mb P$ is in $\cP_\Lambda(X) = 
\cE_\Lambda(T^*X)$, then the vector space $V_\lambda$ from the quiver $R({\mb P}|_{U^2})$
is naturally identified with a stalk of the Morse local system $\cM_\lambda\mb P$.
\end{thm}

The remaining sections of this paper construct the functor $R$ and prove this theorem.
The construction uses the description of the geometry of $\Lambda$ obtained 
in \S2, working inward to describe successively the categories 
$\cE_\Lambda(U^d)$ for $d = 0, 1, 2$. 

\begin{proof}[Proof of main theorems (\ref{Amain} and \ref{Dmain})]
The restriction functor 
\[|_{U^2}\colon\cP_\Lambda(X) = \cE_\Lambda(T^*X) \to \cE_\Lambda(U^2)\] is
a full embedding of categories, by Theorem \ref{GMVcod1}.  Thus 
it is enough to show that every isomorphism
class of objects of $\cE_\Lambda(U^2)$ is in the image of this functor, 
or in other words
that Kashiwara's conjecture holds for our varieties.

We first show that every simple object is in the image.  By Theorem \ref{simples}
an object $(\{V_\lambda\}, \{t_\alpha\}, \{p(\lambda,\lambda')\})$
in $\cA_{k,l}$ or $\cA^s_k$ is simple if and only if 
$\sum_\lambda \dim V_\lambda = 1$.  Thus there is one simple object $S_\mu \in 
\cA_{k,l}$ (respectively $\cA^s_k$) for each $\mu \in \Omega_{k,l}$ (resp. $\Omega^s_{k}$).

We show by induction on $|\mu|$ that $R(\IC(\overline{X_\mu})|_{U^2}) \cong S_\mu$.  If 
$|\mu| = 0$, then $X_\mu$ is a point and the result is clear.  Now, assuming it holds 
for all $\nu$ with $|\nu| < |\mu|$, we need to show that $S = R(\IC(\overline{X_\mu})|_{U^2})$
is irreducible; this is enough, since $V_\mu(S)$ is one-dimensional.  Suppose $S$ were not
irreducible.  Then $S$ has a composition series consisting of one copy of $S_\mu$ and
only using $S_\nu$ for $\nu \le \mu$, and in particular there is a nonzero map between
$S$ and some $S_\nu$ with $|\nu| < |\mu|$.  Using the induction hypothesis, Theorem \ref{abcd},
and Theorem \ref{GMVcod1}, we get a nontrivial morphism between $\IC(\overline{X_\mu})$ and
$\IC(\overline{X_\nu})$, a contradiction.

Finally, to complete the proof, use Theorem \ref{weakKash} and induction to show that
any quiver object is in the image of $R \circ |_{U^2}$.
\end{proof}

\subsection{Cutting by a normal slice}
We prefer to work
on the vector spaces $M_\lambda$, $M^*_\lambda$ and their 
cotangent bundle(s), rather 
than on $X$ directly.  To that end, define
a functor $F_\lambda\colon \cP_\Lambda(X)\to \cP(M_\lambda)$ by
\[F_\lambda{\mb P} = F\epsilon_\lambda^*{\mb P}[-|\lambda|]\]
(note that $\epsilon_\lambda^*{\mb P}$ is conical, so the Fourier transform
can be applied).

Recall the splitting $\Hom(\C^*/W_\lambda, W_\lambda) = M_\lambda \oplus \bar M_\lambda$ 
from \S2.1.
Define a map $\kappa_\lambda\colon T^*M_\lambda \to T^*X$ by sending 
$(x, \xi)$ to $\delta(x,\xi, 0, 0) \in T^*X$, where 
$\delta$ is the inclusion of $T^*M_\lambda \times T^*\bar M_\lambda$ into
$T^*X$ given by differentiating $\tilde\epsilon_\lambda$.

\newcommand{\Lh}{{\hat{\Lambda}}}
Let $\Lh = \Lh_\lambda = \kappa_\lambda^{-1}(\Lambda)$. 
Considered as a subset of $T^*M^*_\lambda = T^*M_\lambda$,
it is the conormal variety to the stratification
of $M^*_\lambda$ by the sets $\epsilon^{-1}(X_{\lambda'})$, $\lambda'\ge\lambda$.
It is contained in the conormal variety to the stratification of $M_\lambda$ by
$B_\lambda$ orbits, but note that it is {\em not\/} in general 
the conormal variety to a stratification of $M_\lambda$.

The next proposition follows easily from 
the fact that $\epsilon_\lambda$
is transverse to the stratification.  Let $\hat{U} \subset T^*M_\lambda$ be an open set, 
and let $U = \delta(\hat{U} \times T^*\bar M_\lambda) \subset T^*X$.

\begin{prop} \label{normal slice} The image of $F_\lambda$ is contained in $\cP_{\Lh}(M_\lambda)$.
More generally, $F_\lambda$ induces an 
equivalence of categories $\cE_\Lambda(U) \to \cE_\Lh(\hat{U})$.
\end{prop} 

\subsection{Codimension zero MPS}
Microlocal perverse sheaves on the smooth part of $\Lambda$ are very easy to describe:

\begin{prop} \label{cod0d} Take $U \subset T^*X$ open with  
$U\cap \Lambda = {\wt{\Lambda}}_\lambda$.  Then $\cE_\Lambda(U)$ 
is equivalent to the category of local systems on ${\wt{\Lambda}}_\lambda$;
one such equivalence is given by the Morse local systems functor 
${\mb P} \mapsto \cM_\lambda{\mb P}$.
\end{prop}
A proof is given in \cite{GMV}.  The Morse local system functor is defined 
in \cite{MV}.

Since ${\wt{\Lambda}}_\lambda \cong {O}_\lambda\times X_\lambda$, 
the Morse local systems of 
${\mb P} \in \cP_\Lambda(X)$ are determined by the following:

\begin{prop} There is a natural isomorphism \[H^{-kl+|\lambda|}(F_\lambda{\mb P})|_{O_\lambda} \cong 
(\kappa_\lambda)^*(\cM_\lambda {\mb P}).\]
\end{prop}

\subsection{MPS through codimension one} \label{MPSc1} We have seen that all the singularities of 
our stratifications are conical (see  Proposition
\ref{hjk}). Thus we can use the techniques of \cite{BG} to describe 
the category $\cE_\Lambda(U^1)$.

Take $\lambda, \lambda'$ with $\lambda\ral^\pi \lambda'$, and 
consider the corresponding codimension one orbit 
$O = O^\pi_\lambda$ in $M_\lambda$.  Embed its normal bundle 
$L =L^\pi_\lambda = T^{}_OM_\lambda$ as a tubular neighborhood
of $O$ which meets only $O$ and the open orbit $O_\lambda$.

Suppose that $\hat{U} \subset T^*M_\lambda$ contains $T^*L$, and
let $U\subset T^*X$ be the corresponding set as described above.
Then define a functor
\[F^\pi_\lambda = F_\lambda|_L\colon
\cE_\Lambda(U) \to \cP_{mon}(L^\pi_\lambda).\] 
Applying Proposition \ref{normal slice} gives immediately

\begin{prop} \label{cod1d} If $\hat U= T^*L$, then $F^\pi_\lambda$ is an equivalence of categories.
\end{prop}

Following Lemma \ref{Pmon}, a perverse sheaf
${\mb P}$ gives rise to local systems \[\cL{\mb P} = F^\pi_\lambda{\mb P}|_{\wt L},\,
\cL'{\mb P} = (FF^\pi_\lambda{\mb P})|_{\wt{L}'}.\] 
In order to combine the local descriptions of microlocal perverse sheaves
from Propositions \ref{cod0d} and \ref{cod1d} 
into a description of $\cE_\Lambda(U^1)$, we need to relate 
$\cL{\mb P}$ and $\cL'{\mb P}$ 
with the Morse local systems $\cM_\lambda{\mb P}$ and $\cM_{\lambda'}{\mb P}$.
The result obtained in \cite{BG} is summarized by the following theorem.

Let $U^\pi_\lambda=(T^*X \setminus \Lambda) \cup \cO_\lambda \cup 
\cO_{\lambda'} \cup \cO^\pi_\lambda$
 be the union of all orbits of 
$\Lambda$ containing $\Lambda^\pi_\lambda$ in their closures.
Let $\Lh^\pi_\lambda = O^\pi_\lambda \times \{0\} \subset \Lh$ be
the codimension one orbit of $\Lh$ lying over $O_\lambda^\pi$. 

We can embed $L$ and $L'$ as tubular neighborhoods of $\Lh^\pi_\lambda$ in
$\kappa_\lambda^{-1}(\cO_\lambda\cup\cO^\pi_\lambda)$ and 
$\kappa_\lambda^{-1}(\cO_\lambda\cup\cO^\pi_{\lambda'})$, respectively.
Define inclusions $i\colon \wt{L} \to {\cO}_\lambda$, 
$i'\colon \wt{L}' \to {\cO}_{\lambda'}$ by restricting these embeddings
and then following them by $\kappa_\lambda$.

\begin{thm} \label{mlpres} \cite{BG} 
There is a one-dimensional local system $\cL_{tw}$ on $\wt{L}'$ 
so that the category $\cE_\Lambda(U^\pi_\lambda)$ is equivalent to the
category whose objects are triples $({\mb P}^\pi_\lambda, \cM, \cM')$, where
${\mb P}^\pi_\lambda\in \cP_{mon}(L^\pi_\lambda)$, and 
$\cM$ and $\cM'$ are local systems on
${\cO}_\lambda$ and ${\cO}_{\lambda'}$ respectively, together
with isomorphisms
\begin{eqnarray*}
\cL({\mb P}^\pi_\lambda) & \cong & i^*(\cM), \\
\cL'({\mb P}^\pi_\lambda) & \cong & (i')^*(\cM') \otimes \cL_{tw},
\end{eqnarray*}
and morphisms are triples of morphisms compatible with these isomorphisms.
Under this equivalence of categories, 
a perverse sheaf ${\mb P}$ is sent to the triple
$(F^\pi_\lambda {\mb P}, \cM_\lambda{\mb P}, \cM_{\lambda'}{\mb P})$. 
\end{thm}

For the spaces we are considering, we have a further simplification:
\begin{prop} The local system $\cL_{tw}$ is trivial.
\end{prop}
\begin{proof} Let $({\mb P}^{\pi}_\lambda, \cM, \cM')$ 
be the triple corresponding to the simple perverse sheaf 
${\mb P} = \IC(\overline{X_{\lambda'}})$.  The Morse
local system $\cM'$ is trivial, so we need to show that 
$\cL'({\mb P^\pi_\lambda})$
is also trivial.  

There is an action of the fundamental group $\pi_1(T)$ on 
${\mb P}$.  This action must be trivial, since ${\mb P}$ is simple, and 
the local system ${\mb P}|_{X_{\lambda'}}$ is trivial.  
The torus $T$ acts on $L$ and $L'$ in a manner compatible with the 
action on the base $O$.  The action of $\pi_1(T)$ on the stalks of 
$\cL'({\mb P^\pi_\lambda})$ must therefore be trivial.

Since $\pi_1(T)$ generates $\pi_1(O)$, it is enough to show that 
$\cL_{tw}$
has trivial monodromy around a loop in a fiber of the
projection $\wt{L}' \to O$.  But by Lemma 5.10 of \cite{BG},
this monodromy is $(-1)^{d-1}$, where $d = 
\dim_\C X_{\lambda'} - \dim_\C X_\lambda$.  The proposition thus follows
from the following result, which is a direct consequence of the 
definition of the relation $\lra$ from \S1.2.
\end{proof}

\begin{lemma} In both $\Omega$ and $\Omega^s$, 
if $\lambda \lra \lambda'$, then
$d$ is odd.
\end{lemma}
\begin{rem} This lemma could also be deduced from \cite{BFL} 
(in type A) and \cite{BF} (in both types).  If it failed to hold
for some $\lambda \ra \lambda'$, then the result of \cite{B} would 
show that
$SS(\IC(\overline{X_{\lambda'}}))$ contains both 
$\overline{T^*_{X_{\lambda'}}X}$ and $\overline{T^*_{X_\lambda}X}$,
contradicting the fact that it is irreducible.  
\end{rem}

Let $A^1$ be the quiver algebra with the same generators and relations 
as $A_{k,l}$ (or $A^s_k$ in the case $X = X^s_k$), except without the
relation 4, and let $\cA^1$ be the
category of finite dimensional 
representations of $A^1$.

\begin{prop} \label{cod1desc} There is an equivalence of categories 
$R\colon \cE_\Lambda(U^1) \to \cA^1$.
\end{prop}
\begin{proof} Essentially this is just gluing together the presentations
of the categories $\cE_\Lambda(U^\pi_\lambda)$ given
by Theorem \ref{mlpres} along the Morse local systems
of Proposition \ref{cod0d}.  

Given ${\mb P} \in \cP_\Lambda(X)$, 
The vector space $V_\lambda$ in the associated quiver 
object will be the stalk of the Morse local system $\cM_\lambda{\mb P}$
at a point $y_\lambda \in {\cO}_\lambda$.  The map $t_\alpha$ is the
action of the monodromy around the loop $e^{i\theta}\mapsto 
\gamma_\alpha(e^{i\theta})\cdot y_\lambda$.  Since these loops generate
$\pi_1({\wt{\Lambda}}_\lambda)$, these maps completely describe
$\cM_\lambda{\mb P}$.  Relations 
1 and 2 (resp. 2s) in the quiver descriptions of $\cA_{k,l}$ and
$\cA^s_k$ give the relations among these loops, by Propositions
\ref{cod0d}, \ref{lsq} and \ref{Dlsq}.

The categories $\cP_{mon}(L^\pi_\lambda)$
have quiver descriptions using Lemma \ref{Pmon}.  We choose the orientation
of $L^\pi_\lambda$ to be $\eta(\beta) = (-1)^{\beta + \half}$ 
times the standard (holomorphic) one, 
where $\pi = (\alpha, \beta)$.  This choice will make
the description of the codimension two relations simpler -- see 
\S\ref{orient}.  
Lemma \ref{Pmon} gives maps between the stalks of the local systems
$\cL({\mb P}^\pi_\lambda)$ and $\cL'({\mb P}^\pi_\lambda)$
at points $z$, $z'$ in $L^\pi_\lambda$ and $(L^\pi_\lambda)'$.  
We can identify these stalks 
with $V_\lambda$ and $V_{\lambda'}$, by fixing
paths from $i(z)$ to $y_\lambda$ in $\cO_\lambda$ and from $i'(z')$ to $y_{\lambda'}$
in $\cO_{\lambda'}$,   
and using Theorem \ref{mlpres}.
This defines the maps  
$p(\lambda, \lambda')$ and $p(\lambda', \lambda)$ in the quiver.

The fact that the maps in 
Lemma \ref{Pmon} are maps of local systems shows
that the maps $t_\alpha$ commute with the maps $p(\lambda, \lambda')$.
The quiver relation 3 (resp. 3s) comes 
from the relations in Lemma \ref{Pmon}, remembering the calculation
of a loop around a point of $O^\pi_\lambda$ in sections 
\ref{pi1} and \ref{Dpi1}.

This defines the required functor $\cE_\Lambda(U^1) \to \cA^1$, and in
fact shows that it factors through an equivalence of categories 
$\cE_\Lambda(\cU) \to \cA^1$, where $\cU$ is the cover of $U^1$
consisting of $U^0$ and the sets $U^\pi_\lambda$ defined in \S4.2.  
The proposition now follows since $\cE$ is a stack.
\end{proof}

\subsection{Codimension two relations}  \label{c2rel}
We now show that the relations 
coming from codimension two strata in $\Lambda$ give exactly the
relation (4) from the quiver categories $\cA_{k,l}$ and $\cA^s_k$.

Take $O = O_A \subset M_\lambda$ a codimension two $B_\lambda$-orbit.
Let $O_1$, $O_2$ be the codimension one orbits containing $O$ in 
their closure (we will assume that there are two; if there is only
one the same discussion works with minor modifications).
Let $\lambda_1, \lambda_2, \lambda'$
be the corresponding elements of $\Omega$, so that $O$, $O_1$, and $O_2$ are
dual to the open orbits in $\epsilon^{-1}(X_{\lambda'})$,
 $\epsilon^{-1}(X_{\lambda_1})$, and $\epsilon^{-1}(X_{\lambda_2})$,
respectively. The 4-tuple 
$(\lambda, \lambda_1, \lambda_2, \lambda')$ forms a diamond, with
$\lambda \ral^{\pi_i} \lambda_i$, $i=1, 2$ --- see Proposition \ref{cod2prop}.

Given an object ${\mb P}$ in $\cE_\Lambda(U^2)$, let ${\mb P}_O$ be the 
restriction 
of $\nu_O^{}F_\lambda{\mb P}[2+|\lambda| - kl]$ to the fiber 
$M_\lambda/T_AO \cong \C^\SO \cong \C^2$,
where $\SO = \Rb \setminus \tau(A)$.  It is a perverse sheaf in
$\cP_{\Lambda_{nc}}(\C^\SO)$,  where $\Lambda_{nc}$ is the conormal 
variety to the normal crossings stratification,
using the decomposition $\C^\SO = L_1 \oplus L_2$, $L_i = T_A\overline{O_i}/T_AO$.  
We also get an induced functor 
$\cE_\Lambda(U^1) \to \cE_{\Lambda_{nc}}(U^1_E)$, where $U^1_E$ is the
set defined in \S\ref{ncpss} which contains all of $\Lambda_{nc}$ 
except the codimension two point.

\begin{prop} \label{cod2desc} In terms of the quiver descriptions of Proposition \ref{cod1desc},
the quiver that Proposition \ref{ncps} associates to ${\mb P}_O$ 
for ${\mb P}$ in $\cE_\Lambda(U^i)$, $i= 1, 2$ is
\[\xymatrix{
V_{\lambda'} \ar@<.5ex>[r] \ar@<.5ex>[d] & 
V_{\lambda_1} \ar@<.5ex>[l]\ar@<.5ex>[d]\\
V_{\lambda_2} \ar@<.5ex>[u] \ar@<.5ex>[r] & 
V_{\lambda} \ar@<.5ex>[l]\ar@<.5ex>[u]}\]
where the maps are the $p(\cdot, \cdot)$ from the algebra $\cA^1$.
\end{prop}

We prove this result in the remaining sections of the paper.

\begin{proof}[Proof of Theorem \ref{abcd}] By Lemma \ref{onec2pt}, 
an object ${\mb P} \in \cE_\Lambda(U^1)$ extends
to an object in $\cE_\Lambda(U^2)$ if and only if each ${\mb P}_O$ 
extends to an object in $\cP_{\Lambda_{nc}}(\C^\SO)$; now use Proposition \ref{ncps}.
\end{proof}

It is clear from Theorem 
\ref{cod1desc} that $V_\lambda$, $V_{\lambda_1}$ and
$V_{\lambda_2}$ and the maps between them 
appear as they do in Proposition \ref{cod2desc}.  What remains is to explain why 
$V_{\lambda'}$ and the four maps involving it appear.

\subsection{Compatibility of orientations} \label{orient}
Our choice of a normal orientation of $M_\lambda$ 
along each codimension one stratum $O^\pi_\lambda$ produces
normal orientations of each irreducible component $\cO_\lambda \subset \Lambda$
along the codimension two orbits contained in it.  In 
order to apply Proposition \ref{ncps} as we do in the statement
of Proposition \ref{cod2desc}, we need to show that
our choice of orientations (\S\ref{MPSc1}) is consistent in 
the sense described after Proposition \ref{ncps}.

We will show this for type A, and leave the type D case
as an exercise. Say that $\pi_i = (\alpha_i, \beta_i)$, and let 
$\pi'_i = (\alpha'_i, \beta'_i)$ be the parent of $\pi_i$ in $\Pi(\lambda)$.

First consider the case where $\pi_1$ and $\pi_2$ are not parent 
and child in $\Pi(\lambda)$.  We then have
$\lambda_1 \ral^{\pi_2} \lambda'$ and $\lambda_2 \ral^{\pi_1} \lambda'$.
If $\pi'_1 \ne \pi'_2$ then $\pi_1'$ is still the parent of $\pi_1$ in 
$\Pi(\lambda_2)$.  Thus the orientations around
$\cO_\lambda^{\pi_1}$ and $\cO_{\lambda_2}^{\pi_1}$ are 
both given by multiplication by
$(\gamma_{\alpha_1}\gamma_{\beta'_1})^{\eta(\beta_1)}$, and
so are compatible.

If $\pi'_1 = \pi'_2$, 
assume WLOG that $\beta_1 < \alpha_2$.  The parent of 
$\pi_1$ in $\Pi(\lambda_2)$ is $(\alpha'_2, \alpha_2)$.  The orientation
around $\cO^{\pi_1}_{\lambda_2}$ is thus given by multiplication by
$(\gamma_{\alpha_1}\gamma_{\alpha_2})^{\eta(\beta_1)}$. 
Since $\gamma_{\alpha_2} = \gamma^{-1}_{\alpha'_2}$ in 
$\pi_1(\cO_{\lambda_2})$,
and $\gamma^{}_{\beta'_1} = \gamma^{-1}_{\alpha'_2}$ in 
$\pi_1({\cO}_{\lambda})$, these define compatible 
orientations.  The other compatibility is easier to check, since the parent
of $\pi_2$ in $\Pi(\lambda_1)$ is 
always of the form $(\tilde\alpha, \beta'_2)$.

Next suppose that $\pi_2$ is the parent of $\pi_1$, and let 
$\pi_3$ be the parent of $\pi_1$.  As was noted in \S\ref{cod2}
and Lemma \ref{diam}, there are two possibilities for $\lambda'$.
Suppose we have $\lambda_1 \ral^{(\beta_1, \beta_2)} \lambda' = \lambda_r$;
the case where $\lambda_1 \ral^{(\alpha_2, \alpha_1)} \lambda'$
is very similar.
The orientation around $\cO^{(\beta_1, \beta_2)}_{\lambda_1}$ is given by 
$(\gamma^{}_{\beta_1}\gamma^{}_{\beta_3})^{\eta(\beta_2)}$, and the orientation
around $\cO^{\pi_2}_\lambda$ is given by 
$(\gamma^{}_{\alpha_2}\gamma^{}_{\beta_3})^{\eta(\beta_2)}$.  The same loops
are generated by $(\gamma^{-1}_{\beta_2}\gamma^{}_{\beta_3})^{\eta(\beta_2)}$
in both cases.
The orientation around $\cO^{\pi_1}_\lambda$ is given by
$(\gamma_{\alpha_1}\gamma^{}_{\beta_2})^{\eta(\beta_1)}$, and the
orientation around $\cO^{(\alpha_2, \alpha_1)}_{\lambda'}$ is
given by $(\gamma_{\alpha_2}\gamma^{}_{\beta_1})^{\eta(\alpha_1)}$.
Compatibility follows from the fact that $\beta_1 - \alpha_1$ is odd.

\subsection{Relating the functors $F_\lambda$} To prove Proposition \ref{cod2desc},
we need to be able to relate the perverse sheaves $F_\lambda{\mb P}$ for 
different $\lambda$.
Fix $\lambda \in \Omega_{k,l}$, and
take an elementary matrix $E = E_{ij}$ in $M^*_\lambda$.
Define $\lambdap \in \Omega_{k,l}$ so that $\epsilon(E) \in X_{\lambdap}$;
it is the partition described by Lemma \ref{elemmat}.  

Define $S = \{(i, s)\in \Rb \mid s \le j\} \cup
\{(r, j) \in \Rb \mid r \le i\}$; we have $T_EO_E = \C^S$.  If we put
$R' = \Rb \setminus S$, and $Z = \C^{R'}$, we have
$Z = (T_EO_E)^\bot \subset M_\lambda$.  
We have a sequence of natural isomorphisms:
\[Z^* \cong M^*_\lambda/T_EO_E = (T_{O_E}M^*_\lambda)|_E \cong
(T_{X_{\lambdap}}X)|_{\epsilon_\lambda(E)} \cong 
(T_{X_\lambdap}X)|_{W_\lambdap} = M^*_\lambdap;\]
the next to last one comes from the trivialization of the normal
bundle $T_{X_\lambdap}X$ induced by the coordinate chart 
$\tilde\epsilon_\lambdap$.


Let 
$\iota\colon Z \to M_{\lambdap}$ be the dual isomorphism.
It will be useful to have an explicit description of it.
Suppose $(\alpha, \beta) = w^{-1}_\lambda(i,j)$ (the map $w_\lambda$
was defined in \S\ref{ncc}).
\begin{lemma} \label{iota} 
If $(r, s) = w^{}_\lambda(\tilde\alpha, \tilde\beta) \in 
R'$, then $\iota(E_{rs}) = E_{r's'}$, where
\[(r',s') = \begin{cases}
w^{}_{\lambdap}(\tilde\alpha, \alpha)&\text{if $\tilde\beta = \beta$},\\
w^{}_{\lambdap}(\beta, \tilde\beta)&\text{if $\tilde\alpha = \alpha$},\\
w^{}_{\lambdap}(\tilde\alpha, \tilde\beta)&\text{otherwise}.
\end{cases}\]
\end{lemma}

The matrix $E\in M^*_\lambda$ gives rise to a linear function 
$f\colon M_\lambda \to \C$, given by $f(A) = A_{ij}$.  We need the
vanishing cycles functor $\phi_f\colon D^b(M_\lambda) \to D^b(f^{-1}(0))$.

\begin{prop} \label{relFl} The perverse sheaf $\phi_f(F_\lambda{\mb P})$ 
is supported on $Z$.  There is a natural isomorphism
\[\phi_f(F_\lambda{\mb P}) \cong \iota^*F_{\lambdap}{\mb P}.\]
\end{prop}
\begin{proof}
Let $\Wb = f^{-1}(0) = (\C E)^\bot$, so 
$\Wb^* \cong M^*_\lambda/\C E$.  Consider the diagram
\[\xymatrix{D^b(M^*_\lambda) \ar[r]^-{\nu^{}_{\C E}} \ar[d]& 
D^b(\C E \oplus \Wb^*) \ar[rr]^{\tilde{h}^*} \ar[dr] \ar[d] & &
 D^b(\Wb^*)\ar[d]\\
D^b(M_\lambda) \ar[r]^-{\nu^{}_\Wb} & 
 D^b(M_\lambda/\Wb \oplus \Wb) \ar[r] & 
 D^b(\C E \oplus \Wb) \ar[r]^{{h}^*} & D^b(\Wb) }\]

Here the leftmost horizontal maps are the indicated specialization
functors, $h(v) = (E, v)$ and $\tilde{h}(v^*) = (E, v^*)$
for $v\in \Wb$, $v^* \in \Wb^*$, and all the
other functors are the appropriate Fourier transforms.  
There are natural isomorphisms making the left square and the 
triangle commute by \cite{BG}, Propositions 
2.3 and 2.4.  The right quadrilateral commutes by the functoriality of
the Fourier transform (\cite{KS}, Proposition 3.7.13).  

The sheaf
$\phi_f(F_\lambda{\mb P})$ is given by taking $\epsilon^*_\lambda({\mb P})
\in D^b(M^*_\lambda)$ 
and mapping down and across the bottom to $D^b(\Wb)$. 
Following the 
upper path, we get a natural isomorphism
\[\tilde h^*\nu^{}_{\C E}\epsilon^*_\lambda{\mb P} \cong 
q^*(\epsilon_{\lambdap})^*{\mb P},\]
where $q\colon \Wb^* \to Z^* \cong M^*_{\lambdap}$ is the composition
of the natural projection and identification maps.  The result follows
now from \cite{KS}, Proposition 3.7.14.
\end{proof} 

\subsection{Proof of Proposition \ref{cod2desc}, continued}
We will first complete the proof of Proposition \ref{cod2desc} in the
case where the codimension two orbit is $O = O^{\pi_1\pi_2}_\lambda$,
for $\pi_1$, $\pi_2$ not parent and child.  Let $(i_1, j_1) = w(\pi_1)$,
$(i_2, j_2) = w(\pi_2)$.  We will apply the results of the previous 
section to the elementary matrix $E = E_{i_1j_1} \in M^*_\lambda$,
so that $\lambdap = \lambda_1$.

Define $S \subset \Rb$, $Z$, $f$, $\iota$ as in the last section.  
Using Lemma \ref{iota}, we have
 $\iota(A^{\pi_1}_\lambda) = {{A}}_{\lambda_1}$
and $\iota(A^{\pi_1\pi_2}_\lambda) = A^{\pi_2}_{\lambda_1}$. 

We want to describe the perverse sheaf ${\mb P}_f = \phi_f(F_\lambda{\mb P})$ 
near the point $A^{\pi_1\pi_2}_\lambda$.  
It will not generally be constructible with respect to the $B_\lambda$-orbits,
but it will be left invariant by the action of the subgroup 
 $B_{\lambda,E} \subset B_\lambda$ of elements which stabilize 
$[E] \in \PP M^*_\lambda$.
Thus it is
enough to study ${\mb P}_f$ restricted to a normal slice to the orbit
$B_{\lambda,E}A_\lambda^{\pi_1\pi_2}$.  One such slice is given by
\[N = A^{\pi_1\pi_2}_\lambda + \C^{S \cup \{(i_2,j_2)\}}.\]

Using the 
isomorphism ${\mb P}_f \cong \iota^*F_{\lambda_1}{\mb P}$ 
given by Proposition \ref{relFl}, we see that 
on $N \cap Z = \C^{\{(i_2, j_2)\}}$, ${\mb P}_f$ is given by 
the quiver $A_{\lambda_1} \leftrightarrows A_{\lambda'}$, 
where the maps are $p(\lambda_1, \lambda')$ and $p(\lambda', \lambda_1)$.
Thus it will be enough give a natural isomorphism (defining $S_O$, ${\mb P}_O$ as in 
\S\ref{c2rel})
\[{\mb P}_f|_N[i_1+j_1 + |\lambda| - kl] \cong \phi_h{\mb P}_O,\]
where $h\colon \C^\SO \to \C$ is given by $h(A) = A_{i_1j_1}$ and
we identify $N \cap Z$ with $\C^{\{(i_2,j_2)\}}$ in the obvious way.

\begin{figure}
\begin{center}
\leavevmode
\hbox{
\epsfxsize=1.8in
\epsffile{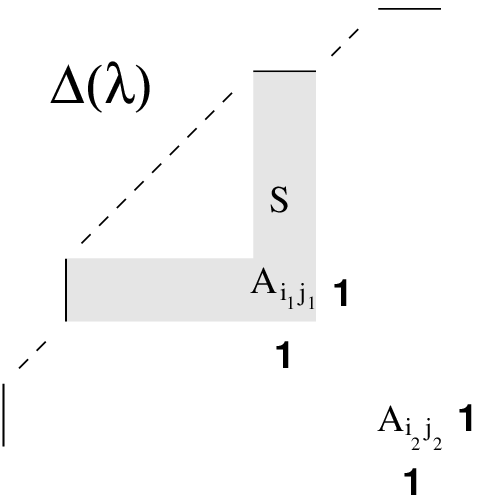}}
\end{center}
\caption{}
\label{lastfig}
\end{figure}

By Proposition
\ref{Mlorbs}, the divisor $Y_m = \overline{O^{\pi_m}_\lambda}$, $m = 1, 2$, 
is given by the vanishing of the determinant
of the submatrix with lower right corner at $(i_m, j_m)$
and upper left corner at $w_{\lambda_m}(\pi_m)$.  
Assume for the moment that 
$\pi_2 \not< \pi_1$. 
A simple calculation gives
\[Y_1 \cap N = \{A \in N \mid A_{i_1j_1} + Q(A_{S^-}) = 0\},\,\text{and}\] 
\[Y_2 \cap N = \{A \in N \mid A_{i_2j_2} = 0\}.\]
Here $Q$ is the nondegenerate quadratic form on $\C^{S^-}$, 
$S^- = S \setminus \{(i_1,j_1)\}$ given by the formula
$Q = \sum A_{ij_1}A_{i_1j}$, where the sum is over all $(i,j) \in\Supp 
A^{\pi_1\pi_2}_\lambda$ with $(i,j) < (i_1,j_1)$
(see Figure \ref{lastfig}).

Applying the change of coordinates 
\begin{align*}
A'_{i_1j_1} & = A_{i_1j_1} + Q(A_{S^-}), \\
A'_{rs} & = A_{rs}, \;\, (r,s) \ne (i_1,j_1)
\end{align*}
gives $Y_m \cap N = \{A' \in N \mid A'_{i_mj_m} = 0\}$ and 
$f(A') = A'_{i_1j_1} - Q(A'_{S^-})$.  The claim follows
now from the following lemma, which generalizes a result \cite{ST} 
of Sebastiani and Thom;
see \cite{M} for a proof.

\begin{lemma} \label{vcstab} Take a polynomial $f \colon \C^n \to \C$,
and let $\tilde{f}\colon \C^n \times \C^k \to \C$ be given by 
\[\tilde{f}({\mb v}; x_1, \dots, x_k) = f({\mb v}) + x_1^2 + \dots + x_k^2.\]
Let $j\colon f^{-1}(0) \to \tilde{f}^{-1}(0)$ be given by 
$j({\mb v}) = ({\mb v}, 0)$, and let $p_1\colon \C^n \times \C^k \to \C^n$
be projection on the first factor.
If ${\mb A} \in D^b(\C^n)$, 
there is a natural isomorphism
\[\phi_{\tilde{f}}p^*_1{\mb A}[k] \cong Rj_*\phi_f{\mb A}.\]
\end{lemma}

Proposition \ref{cod2desc} now follows; just switch $\pi_1$ and $\pi_2$
and repeat the argument.
If $\pi_2 < \pi_1$, a similar argument works, but the equation for
$Y_1 \cap N$ has a slightly different form.  We leave the details for the 
reader.  We also omit the type D case; the 
varieties $Y_m$ have equations given by the vanishing
of the appropriate Pfaffians, but otherwise the argument is 
essentially the same.

\subsection{Proof of Proposition \ref{cod2desc}, concluded}
It only remains to handle the cases $O = O^{\pi_1\pi_2,l}_\lambda$
and $O = O^{\pi_1\pi_2,r}_\lambda$,
where $\pi_1, \pi_2 \in \Pi(\lambda)$, and $\pi_2$ is the parent of
$\pi_1$.  The argument is similar for 
the two cases; we will give it for the first one.  Assume WLOG that
$\pi_2$ is the parent of $\pi_1$, and put $(i_m, j_m) = w_\lambda(\pi_m)$.
Apply Proposition \ref{relFl},
using the elementary matrix $E = E_{i_2j_1}$; we get  
$\lambdap = \lambda'$.  Also
$\iota(A^{\pi_1\pi_2,l}_\lambda) = A^{\pi'}_\lambda$, 
where $\pi' = (\beta_1, \beta_2)$ satisfies $\lambda' \ral^{\pi'} \lambda_2$.

Put ${\mb P}_f = \phi_f(F_\lambda{\mb P}) \cong \iota^*F_{\lambda'}{\mb P}$.
Defining $S$, $Z$ as before, we again see that to describe ${\mb P}_f$
near $A^{\pi_1\pi_2,l}_\lambda$
 it is enough to describe its restriction to the normal
slice \[N = A^{\pi_1\pi_2,l}_\lambda + \C^{S \cup \{(i_2,j_2)\}}.\]
  Since ${\mb P}_f$ is
supported on $Z = \C^{\Rb \setminus S}$, ${\mb P}_f|_N$ is 
a monodromic perverse sheaf on 
$N\cap Z = A^{\pi_1\pi_2,l}_\lambda + \C^{\{(i_2,j_2)\}} \cong \C.$
We will show that there is a natural isomorphism
\[ {\mb P}_f|_N \cong F\phi_h {\mb P}_O,\]
letting $h(A) = A_{i_1j_1}$, and
identifying $N\cap Z$ with $\C^{\{(i_2,j_2)\}}$ as before.

Again putting $Y_m = \overline{O^{\pi_m}_\lambda}$ for $m = 1,2$, we get 
\[Y_1 \cap N = \{A \in N \mid A_{i_1j_1} = 0\},\;\text{and}\]
\[Y_2 \cap N = \{A \in N \mid A_{i_2j_1} - A_{i_1j_1}A_{i_2j_2} +
 Q(A_{S^-}) = 0\},\]
where $Q$ is a nondegenerate quadratic form on $\C^{S^-}$, 
$S^- = S \setminus \{(i_1,j_1), (i_2, j_2)\}$.

Define a change of variables as follows:
\begin{align*}
A'_{i_2j_1} & = A_{i_2j_1} - A_{i_1j_1}A_{i_2j_2} + Q(A_{S^-}), \\
A'_{rs} & = A_{rs}, \;\, (r,s) \ne (i_2,j_1).
\end{align*} 
The result now follows from Lemma \ref{vcstab} and the following
lemma.

Let $V, W$ be one-dimensional 
complex vector spaces, and take a biconic
sheaf ${\mb A} \in D^b(V \times W)$, i.e. let ${\mb A}$ be constructible 
with respect to the normal crossings stratification.
Fix an element $w^* \in W^*$. Let 
$p_{12} \colon V \times W \times V^* \to V \times W$ be the projection,
let $g\colon V \times W \times V^* \to \C$ be given by $g(v, w, v^*) = 
\langle v, v^* \rangle + \langle w, w^* \rangle$. 
Let $i\colon
V^* \to V^* \times W^*$ and $j \colon V^* \to g^{-1}(0)$ be given by
$i(v^*) = (v^*, w^*)$ and $j(v^*) = (0, 0, v^*)$.

\begin{lemma} 

There is a natural isomorphism
\[j^*\phi_g(p^*_{12}{\mb A})[1] \cong i^*F{\mb A}.\]
\end{lemma}
\begin{proof} 
Consider $X = V\times W \times V^*$ as a vector bundle over the
base $B = V^*$, so that $E = g^{-1}(0)$ is a sub-vector bundle.  Then
the vanishing cycles functor $\phi_g$ is naturally equivalent
to $s^*F\nu_E[-1]$, where $\nu_E\colon X \to T_EX = E \times_B X/E$
is the specialization functor, $F$ is the Fourier transform
$D^b(T_EX) \to D^b(T^*_EX)$, and $s\colon E \to T^*_EX$ is
given by $s(\zeta) = (\zeta, dg_\zeta)$.

Now we do a diagram chase similar to the one in Proposition \ref{relFl}
Consider the following diagram of functors:
\[\xymatrix{
D^b(V\times W) \ar[r]^-{p^*_{12}} & D^b(X) \ar[r]^-{\nu^{}_E} & 
D^b(T_EX) \ar[r] & D^b(T^*_EX) \ar[r]^-{s^*} & D^b(E)\\
D^b(V^* \times W^*) \ar[u] \ar[r]^-{p^*_{12}} & D^b(X^*) \ar[u] 
\ar[r]^-{\nu^{}_{E^\bot}} & D^b(T_{E^\bot}X^*) \ar[u] \ar[ur]\ar[r]^-{(s')^*} &
D^b(E^*)\ar[ur]
}\]
where $E^\bot$ is the subbundle of $X$ which annihilates $E$, 
the map $s'\colon E^* \to T_{E^\bot}X^* = E^* \times_B (X/E)^*$ is
given by $s'(\xi) = (\xi, dg_{(0,v^*,0)})$ (putting 
$v^* = \pi_{E^*}(\xi)$), and all the 
unmarked arrows are the appropriate Fourier transform functors.
Arguing as in Proposition \ref{relFl}, we see that
that there is a natural isomorphism
\[j^*\phi_g(p^*_{12}{\mb A})
 \cong j^*F(s')^*\nu^{}_{E^\bot}p_{12}^*F{\mb A}.\]

Now, using the assumption that the vector bundles $V$ and $W$
are one-dimensional, a simple argument shows
 that there is a natural isomorphism
\[(s')^*\nu^{}_{E^\bot}p_{12}^*F{\mb A} \cong (\pi_{E^*})^*i^*F{\mb A}. \]
But then for any ${\mb C} \in D^b(B)$, $j^*F(\pi_{E^*})^*{\mb C}$ is 
naturally isomorphic to ${\mb C}[-1]$, which
gives the result.
\end{proof}

\end{document}